\documentclass[hidelinks,onefignum,onetabnum]{siamart250211}

\usepackage{bm,bbm}
\usepackage{subfigure}
\usepackage{color}
\usepackage{array, makecell}
\usepackage{enumitem}

\usepackage{lipsum}
\usepackage{amsfonts}
\usepackage{graphicx}
\usepackage{epstopdf}
\usepackage{algorithmic}
\ifpdf
  \DeclareGraphicsExtensions{.eps,.pdf,.png,.jpg}
\else
  \DeclareGraphicsExtensions{.eps}
\fi

\newsiamremark{remark}{Remark}
\newsiamremark{hypothesis}{Hypothesis}
\crefname{hypothesis}{Hypothesis}{Hypotheses}
\newsiamthm{claim}{Claim}
\newsiamremark{fact}{Fact}
\crefname{fact}{Fact}{Facts}

\headers{Dynamical low-rank method for the Wigner equation}{S. Shao and Y. Shao}

\title{Dynamical low-rank methods for the Wigner equation I: separable difference potential\thanks{Submitted to the editors DATE.
\funding{This research was supported by the National Natural Science Foundation of China (Nos. 12325112, 12288101) and the High-performance Computing Platform of Peking University..}}}

\author{Sihong Shao\thanks{CAPT, LMAM and School of Mathematical Sciences, Peking University, Beijing 100871, China
  (\email{sihong@math.pku.edu.cn}).}
\and Yuehan Shao\thanks{CAPT, LMAM and School of Mathematical Sciences, Peking University, Beijing 100871, China 
  (\email{joshao@pku.edu.cn}).}
}

\usepackage{amsopn}

\newcommand\D{\textup{d}}

\def\mi{\mathrm{i}}
\def\me{\mathrm{e}}

\def\ba{\bm{a}}
\def\bb{\bm{b}}
\def\bx{\bm{x}}

\def\bk{\bm{k}}
\def\bp{\bm{p}}

\def\by{\bm{y}}
\def\bb{\bm{b}}

\def\pdo{{\rm \Psi} \textup{DO}}

\crefname{subsection}{Section}{Sections}
\Crefname{subsection}{Section}{Sections}

\hypersetup{
	pdftitle={Dynamical low-rank methods for the Wigner equation I: separable difference potential},
	pdfauthor={S. Shao and Y. Shao}
}

\begin{document}

\maketitle

\begin{abstract}
  Recent advances in dynamical low-rank approximation (DLRA) have demonstrated its effectiveness in high-dimensional simulations. 
  However, existing DLRA algorithms still face significant challenges when handling systems that involve complex collision terms, including the pseudo-differential operator ($\pdo$) in the Wigner equation, a representative operator characterized by nonlocality. 
  It is deserving to carry out a series of works to develop the DLRA algorithms for solving the Wigner equation. 
  As the first step in this series of works, we propose an efficient DLRA algorithm for the Wigner equation, using a separable decomposition of the difference potential. 
  We combine this separable assumption with two often-used truncations of $\pdo$, namely $\mathcal{K}$-truncation and $\mathcal{Y}$-truncation, to obtain a kind of separated representation of $\pdo$. 
  Complexity analysis and several challenging experiments, including harmonic oscillators, Gaussian barrier scattering, electron-electron scattering, and a Helium-like system, all of which satisfy the separable assumption, confirm that the proposed DLRA algorithm has significant advantages, achieving a reduction in computational effort by one to two orders of magnitude in both runtime and memory requirements compared to the full-grid approach. 
  It is worth noting that, even in the absence of a predetermined low-rank structure for the solution, DLRA can still serve as a numerical scheme that balances efficiency and accuracy.
\end{abstract}



\section{Introduction}
Since its introduction in 1932 \cite{wigner1932quantum}, the Wigner equation provides a convenient way to characterize fundamental phase-space formulation of quantum mechanics, analogous to the Boltzmann equation in classical kinetic theory. 
Although the Wigner function may exhibit negative values, this quasi-distribution function remains an important way to understand and forecast quantum phenomena in many fields, such as nanoelectronics \cite{vandeput2017efficient, gehring2005wigner}, high energy density physics \cite{graziani2014kinetic}, quantum optics \cite{schleich2015quantum}, quantum transport \cite{jacoboni2001quantum,jacoboni2004wignerfunction}, etc. 

\subsection{The Wigner equation}
The Wigner equation with a time-independent potential considered in this paper is
\begin{equation}
	\frac{\partial}{\partial t} f(\bx, \bk, t)=-\frac{\hbar \bk}{m} \cdot \nabla_{\bx} f(\bx, \bk, t)+\Theta_V[f](\bx, \bk, t):=\mathcal{L}f, \label{wigner}
\end{equation}
where $\bx \in \Omega_{\bx} \subseteq \mathbb{R}^d$ is the spatial variable, $\bk \in \Omega_{\bk} \subseteq \mathbb{R}^d$ is the Fourier conjugated wave vector, $d$ is the dimension, $m$ is the mass, $\hbar$ is the reduced Planck constant and the nonlocal pseudodifferential operator ($\pdo$) is defined as
\begin{align}
	\Theta_V[f](\bx, \bk, t)&=\int_{\mathbb{R}^d} f\left(\bx, \bk^{\prime}, t\right) V_w(\bx,\bk-\bk^{\prime}) \D \bk^{\prime},\label{Vw to theta}\\
	V_w(\bx,\bk)&=\frac{1}{\mi \hbar(2 \pi)^d}\int_{\mathbb{R}^d} \me^{-\mi \bk\cdot\by} D_V(\bx, \by)\D \by, \label{Vw}\\
	D_V(\bx, \by) &= V\left(\bx + \frac{\by}{2}\right) - V\left(\bx - \frac{\by}{2}\right).
\end{align}
Here $D_V(\bx, \by)$ is central difference of potential function $V(\bx)$, the Wigner kernel $V_w(\bx,\bk)$ is Fourier transformation of $D_V(\bx, \by)$ and $\Theta_V[f](\bx, \bk, t)$ is the convolution of $f(\bx, \bk, t)$ and $V_w(\bx,\bk)$. 
We can obtain another form of $\pdo$ by changing the order of integration
\begin{align}
	\Theta_V[f](\bx, \bk, t)&=\frac{1}{\mi\hbar} \int_{\mathbb{R}^d}\me^{-\mi \bk\cdot\by} D_V(\bx, \by) \rho(\bx,\by,t) \D \by, \label{rho to theta}\\
	\rho(\bx,\by,t)&=\frac{1}{(2\pi)^d}\int_{\mathbb{R}^d} f\left(\bx, \bk^{\prime}, t\right)  \me^{\mi\bk^{\prime}\cdot\by}\D \bk^{\prime}. \label{rho}
\end{align}
The $\pdo$ is the most complex and important component in the Wigner equation, representing the fundamental distinction between the classical Vlasov equation and the Wigner equation.
This equation evolves in a $2d$-dimensional phase space, thus confronting the severe challenge of the curse of dimensionality (CoD) along with the nonlocality. 

Despite the long history and broad application of the Wigner equation, its numerical solution continues to present significant challenges nowadays. 
Over the past three decades, significant progress has been made in developing Wigner equation solvers, which can be roughly classified into particle-based stochastic and grid-based deterministic approaches.
The particle-based stochastic approaches such as Monte Carlo method \cite{kosina2003monte, sellier2015simulation} and stochastic particle method \cite{shao2019branching} typically achieve only $\mathcal{O}(N^{-1/2})$ convergence due to their inherent statistical nature, where $N$ is the number of samples. 
In contrast, the grid-based deterministic approaches exhibit superior convergence rates because of their concise mathematical principle, including the finite difference method \cite{frensley1989effect}, operator splitting method \cite{ringhofer1990spectral,arnold1995operator}, spectral element method (SEM) \cite{shao2011adaptive,xiong2016advectivespectralmixed,chen2022higherorder,xiong2023characteristicspectralmixed}, 
moment method \cite{li2014numerical}, weighted essentially non-oscillatory (WENO) method \cite{dorda2015wenosolver}, etc. 
While CoD confined deterministic methods to low-dimensional problems, rapid advances in high-performance computing, alongside techniques such as tensor approximations that mitigate CoD, have enabled deterministic methods to tackle 5D and 6D problems. 
The dynamical low-rank approximation (DLRA) method we consider next serves as one such approach for mitigating CoD.

\subsection{DLRA: the projector-splitting integrator}
Emerging from the Dirac-Frenkel time-dependent variational principle \cite{dirac1930note}, the DLRA method \cite{koch2007dynamical} has evolved into a powerful class of numerical integrators. 
This method has been widely used in kinetic theory, including the Vlasov equation \cite{einkemmer2018lowrank,einkemmer2021mass}, Boltzmann equation \cite{einkemmer2021efficient,hu2022adaptive} and radiation transfer equation \cite{peng2020lowrank,einkemmer2021asymptoticpreserving}. 
The main idea is to constrain the unknown function $f$ to a low-rank manifold, then deriving evolution equations of basis functions through tangent space projection. 
This approach reduces the original $2d$-dimensional problem to several $d$-dimensional subproblems, effectively alleviating CoD in traditional grid-based methods.

The DLRA method begins by assuming that the Wigner function $f(\bx, \bk, t)$ admits a separable low-rank representation of the following form
\begin{equation}
	f(\bx, \bk, t)\approx \sum_{i,j=1}^r X_i(\bx,t)K_j(\bk,t)S_{ij}(t), \label{low-rank}
\end{equation}
where $r$ is the rank, the coefficient matrix $S=(S_{ij})$ is invertible, the basis function $X_i(\bx,t) \in L^2(\Omega_{\bx})$ and $K_j(\bk,t) \in L^2(\Omega_{\bk})$ are orthonormal
\begin{align}
	\left\langle X_i , X_j \right\rangle_{\bx} = \delta_{ij}, \quad  \left\langle K_i , K_j \right\rangle_{\bk} = \delta_{ij}, 
\end{align}
with $\left\langle \cdot , \cdot  \right\rangle_{\bx}$ and  $\left\langle \cdot , \cdot  \right\rangle_{\bk}$ being the standard $L^2$ inner products on $L^2(\Omega_{\bx})$ and $L^2(\Omega_{\bk})$, respectively, and $\delta_{ij}$ denotes the Kronecker delta. 
Therefore, for fixed time $t$, the Wigner function $f(\bx, \bk, t)$ lies in the fixed rank-$r$ manifold
\begin{equation}
	\begin{aligned}
		\mathcal{M}_r=\biggl\{g\in L^{2}(\Omega)\mid g(\bx, \bk)=\sum_{i, j=1}^r X_{i}(\bx) S_{i j} K_{j}(\bk) \text{ with invertible } S=\left(S_{i j}\right) , \\
		X_{i} \in L^{2}\left(\Omega_{\bx}\right), K_{j} \in L^{2}\left(\Omega_{\bk}\right) \text{ with } \left\langle X_{i}, X_{j}\right\rangle_{\bx}=\delta_{i j},\left\langle K_{i}, K_{j}\right\rangle_{\bk}=\delta_{i j}\biggr\}.
	\end{aligned}\label{Mr}
\end{equation}
Then project $\mathcal{L}f$ in \cref{wigner} onto the tangent space of $\mathcal{M}_r$, yielding the reduced system 
\begin{align}
	\frac{\partial}{\partial t} f(\bx, \bk, t)=\mathcal{P}_f(\mathcal{L}f), \label{reduced Wigner}
\end{align}
where the projection operator $\mathcal{P}_f$ acting on a generic function $g$ can be written as \cite{koch2007dynamical}
\begin{align}
	\mathcal{P}_f g= \mathcal{P}_{\bar{K}}g-\mathcal{P}_{\bar{K}}\mathcal{P}_{\bar{X}}g+\mathcal{P}_{\bar{X}}g. \label{PfgPxPk}
\end{align}
Here, $\bar{K}$ and $\bar{X}$ denote the vector spaces spanned by the basis functions $K_j$ and $X_i$,
\begin{align}
	\bar{K}=\operatorname{span}\{K_j\mid j=1,2,\cdots,r\},	\qquad \bar{X}=\operatorname{span}\{X_i\mid i=1,2,\cdots,r\},
\end{align}
and $\mathcal{P}_{\bar{K}}, \mathcal{P}_{\bar{X}}$ are the orthogonal projectors onto $\bar{X}$ and $\bar{K}$, respectively,
\begin{align}
	\mathcal{P}_{\bar{K}}g= \sum_{j=1}^{r} \left\langle K_j, g\right\rangle_{\bk} K_j, \qquad \mathcal{P}_{\bar{X}}g=\sum_{i=1}^r X_i\left\langle X_i, g\right\rangle_{\bx}.
\end{align}

To avoid the possible ill-conditioning of the coefficient matrix $S$ \cite{koch2007dynamical,kieri2016discretized}, the projector-splitting integrator \cite{lubich2014projectorsplitting} substitutes \eqref{PfgPxPk} into the reduced Wigner system \cref{reduced Wigner} and decomposes \cref{reduced Wigner} into three subproblems: 
\begin{align}
	\frac{\partial}{\partial t} f(\bx, \bk, t) & =\mathcal{P}_{\bar{K}}\left(\mathcal{L}f\right),\label{PK}                          \\
	\frac{\partial}{\partial t} f(\bx, \bk, t) & =-\mathcal{P}_{\bar{K}}\mathcal{P}_{\bar{X}}\left(\mathcal{L}f\right), \label{PKPX} \\
	\frac{\partial}{\partial t} f(\bx, \bk, t) & =\mathcal{P}_{\bar{X}}\left(\mathcal{L}f\right).\label{PX}
\end{align}
The decomposition of the projection operator $\mathcal{P}_f$ forms the basis of projector-splitting integrator. 
Exploiting orthogonality conditions, we can further simplify each subproblems and proceed in the following three substeps, respectively:
\subsection*{$Z$-step} We can represent the low-rank assumption for $f$ in \cref{low-rank} as $f(\bx, \bk, t)=\sum_{j=1}^r Z_j(\bx,t)K_j(\bk,t)$ with $ Z_j(\bx,t)=\sum_{i=1}^r X_i(\bx,t) S_{ij}(t),$ then rewrite Eq.\cref{PK} as
\begin{align*}
	\sum_{j=1}^r\left(\frac{\partial}{\partial t}Z_jK_j+Z_j\frac{\partial}{\partial t}K_j\right)=\sum_{j=1}^r\left\langle K_{j},\mathcal{L}f\right\rangle_{\bk} K_{j}.
\end{align*}
Making an important observation that $\{K_j(\bk, t) \}_{1\le j \le r}$ are independent of time $t$ (see \cite[Lemma 3.1]{lubich2014projectorsplitting}) in this substep, and using orthogonality conditions $\left\langle K_i , K_j \right\rangle_{\bk} = \delta_{ij}$, we derive the following differential equation for $Z_j(\bx,t)$ as
\begin{align}
	\frac{\partial}{\partial t}Z_j(\bx,t)= \left\langle K_{j},\mathcal{L}f\right\rangle_{\bk},\qquad j=1,...,r. \label{Z-step}
\end{align}

\subsection*{$S$-step} We can argue similarly to obtain that Eq.\cref{PKPX} is equivalent to
\begin{align}
	\frac{\partial}{\partial t} S_{ij}(t)=\left\langle X_{i} K_{j}, \mathcal{L}f \right\rangle_{\bx,\bk},\qquad i,j=1,...,r. \label{S-step}
\end{align}
During this substep, both $\{K_j(\bk, t) \}_{1\le j \le r}$ and $\{X_i(\bx, t) \}_{1\le i \le r}$ remain unchanged.

\subsection*{$L$-step} Similarly, we can represent \cref{low-rank} as $f(\bx, \bk, t)=\sum_{i=1}^r X_i(\bx, t)L_i(\bk, t)$ with $L_i(\bk, t)=\sum_{j=1}^rS_{ij}(t)K_j(\bk, t).$ By similar arguments, Eq.\cref{PX} is equivalent to
\begin{align}
	\frac{\partial}{\partial t} L_i(\bk, t)=\left\langle X_{i} , \mathcal{L}f \right\rangle_{\bx},\qquad i=1,...,r. \label{L-step}
\end{align}
During this substep, $\{X_i(\bx, t) \}_{1\le i \le r}$ remains unchanged.

This projector-splitting integrator has been shown to be robust to small singular values of $S$ and own first order rate \cite{kieri2016discretized}. 
For more details, we refer the reader to \cite{einkemmer2018lowrank} which provides similar approach for the Vlasov-Poisson equation. 

\subsection{The main challenge}
\label{The main challenge}
While \cref{PK}-\cref{PX} appear to offer a framework for applying DLRA, the efficiency relies on the separability of the right-hand side $\mathcal{L} f$.
The vast majority of problems solved using DLRA possess separability in their right-hand side terms, such as the Vlasov equation \cite{einkemmer2018lowrank,einkemmer2021mass} and the radiative transfer equation \cite{peng2020lowrank,einkemmer2021asymptoticpreserving}. 
Some studies utilize the steady-state properties to simplify the problem, for instance, the Maxwellian distribution is employed to reduce complexity when solving the Boltzmann-BGK equation \cite{einkemmer2021efficient}.
However, this approach is not well suited for solving non-steady-state problems.

The so-called separability here refers to the fact that, under the low-rank assumption \cref{low-rank} satisfied by $f$, $\mathcal{L} f$ can also be naturally expressed in a form that separates the variables $\bx$ and $\bk$. 
Separability allows the projection to be computed directly via the low-rank components $X_i,S_{ij},K_j$, without the full tensor $f$. 
For example, the advection term $\frac{\hbar \bk}{m} \cdot \nabla_{\bx} f$ possesses separability, since 
\begin{align*}
	\frac{\hbar \bk}{m} \cdot \nabla_{\bx} f(\bx, \bk, t)=\frac{\hbar}{m}\sum_{i',j'=1}^r  \nabla_{\bx} X_{i'}(\bx,t) \cdot (\bk K_{j'}(\bk,t)) S_{i'j'}(t),
\end{align*}
and then the projection onto the $\{K_j(\bk, t) \}_{1\le j \le r}$  becomes 
\begin{align*}
	\left\langle K_{j},\frac{\hbar \bk}{m} \cdot \nabla_{\bx} f(\bx, \bk, t)\right\rangle_{\bk}= \frac{\hbar}{m} \sum_{i',j'=1}^r\nabla_{\bx} X_{i'}(\bx,t) \cdot \left\langle K_{j}, \bk K_{j'}(\bk,t) \right\rangle_{\bk}S_{i'j'}(t),
\end{align*}
where $\bx$ and $\bk$ dependent terms can be computed separately.
In contrast, the intrinsic nonlocality of $\pdo$ introduces strong coupling between the variables $\bx$ and $\bk$, preventing it from possessing the separability.
Consequently, evaluating its projection would typically necessitate reconstructing the full tensor, which undermines the efficiency aims of the low‑rank approach.
This strong coupling within $\pdo$ represents the main difficulty to applying DLRA to the Wigner equation.

\subsection{Organization}
The rest of this paper is organized as follows. 
\Cref{Separated approximation for pdo} introduces a separable assumption for the difference potential $D_V(\bx, \by)$ as shown in \cref{DV low-rank}. 
By introducing this separable assumption, $\pdo$ can be transformed into a separated form through $\mathcal{K}$-truncation and $\mathcal{Y}$-truncation.
The strong coupling between $\bx$ and $\bk$ in $\pdo$ can be transferred to the separation rank $R$. 
Regarding the applicability of this assumption, it effectively covers scenarios where the potential $V(\bx)$ is real analytic, including typical potentials such as the harmonic oscillator, Gaussian barrier, and soft Coulomb discussed in subsequent numerical examples.

In \Cref{DLRA for wigner}, we derive the evolution equations for basis functions from the Wigner equation, detail their numerical implementation and analyze the overall complexity of the algorithm. 
In \Cref{Numerical experiments}, we validate the efficiency of our algorithm through a series of typical numerical examples. 
We first choose several non-stationary and atypical low-rank numerical experiments such as harmonic oscillator and Gaussian barrier scattering. 
Our numerical experiments indicate that the DLRA method can still yield numerical solutions of sufficient accuracy, even when the solution lacks an obvious low‑rank structure.
We further validate the method through more complex numerical examples, including electron-electron scattering and Helium-like systems, successfully capturing fundamental quantum phenomena such as the Pauli exclusion principle and the Heisenberg uncertainty principle \cite{heisenberg1927uber}. 
The numerical results confirm that the proposed algorithm can achieve a reduction in computational effort by one to two orders of magnitude in both runtime and memory requirements compared to the full-grid approach.
These results clearly demonstrate the broad applicability and promising potential of the DLRA method in quantum dynamical simulations. 
\Cref{Conclusion and discussion} makes a summary and discusses future research directions.

\section{Separable assumption}
\label{Separated approximation for pdo}


The intrinsic nonlocality of $\pdo$ in \eqref{Vw to theta} places the first obstacle for any numerical solver and some kinds of truncation have been introduced in the literature.  
If the potential $V(\bx)$ admits a Taylor expansion in $\bx$-space, then the difference potential $D_V(\bx, \by)$ can be written as 
\begin{align}
	D_V(\bx, \by)=\sum_{l=0}^{+\infty}\frac{2}{(2l+1)!}\nabla_{\bx}^{(2l+1)}V\left(\bx\right) \left(\frac{\by}{2}\right)^{2l+1}.\label{Dv_Taylor_expansion}
\end{align}
Substituting \cref{Dv_Taylor_expansion} into \cref{Vw} and using the basic properties of Fourier transform yields direatly the so-called Moyal expansion of $\pdo$
\begin{align}
	\Theta_V^{M}[f](\bx, \bk, t)=\sum_{l=0}^{+\infty} \frac{(-1)^{l}}{(2l+1)!}\left(\frac{\hbar}{2}\right)^{2l}\nabla_{\bx}^{2l+1}V\cdot\nabla_{\bp}^{2l+1}f, \label{MT}
\end{align}
where $\bp=\hbar\bk$ is the momentum. 
In particular, the first term in the summation on the right-hand side is $\nabla_{\bx}V\cdot\nabla_{\bp}f$, indicating that the Wigner equation reduces to the classical Vlasov equation when $\hbar$ tends to zero. In the early stage, one has to truncate the infinite series summation \cref{MT} to a finite number of terms \cite{HugMenkeSchleich1998I,HugMenkeSchleich1998II}. Obviously, the resulting scheme possesses certain separability and the primary reason is that in Eq. \cref{Dv_Taylor_expansion} for the difference potential $D_V(\bx,\by)$, the variables $\bx$ and $\by$ are separated. 
Based on this observation, we introduce the separable assumption for the difference potential $D_V(\bx, \by)$ 
\begin{align}
	D_V(\bx, \by) \approx \sum_{s=1}^{R} D_{X,s}(\bx) D_{Y,s}(\by), \label{DV low-rank}
\end{align}
where $R$ denotes the separation rank, and $D_{X,s},D_{Y,s}$ are functions that depend only on $\bx,\by$, respectively. 

Nowadays, the following two kinds of truncation, named after the truncation of the outermost infinite integral in $\pdo$, are often used \cite{chen2022higherorder,chen2019high,chen2019numerical,xiong2016advectivespectralmixed}:
\begin{itemize}[leftmargin=1.5em, itemindent=0em, labelsep=0.3em, itemsep=0.3em, topsep=0.3em]
\item $\mathcal{Y}$-truncation: starting from \cref{rho to theta} and \cref{rho}, computations are performed via the intermediate variable $\rho$;
\item $\mathcal{K}$-truncation: starting from \cref{Vw to theta} and \cref{Vw}, computations are performed via the intermediate variable $V_w$;
\end{itemize}
After applying the separable difference potential assumption \cref{DV low-rank}, we have a separated approximation for $\pdo$
\begin{align}
	\Theta_V[f](\bx, \bk, t)&=\sum_{s=1}^{R}\left[ \sum_{i,j=1}^r \overline{X}_{i,s}(\bx,t)S_{ij}(t)\overline{K}_{j,s}(\bk,t) \right] \label{DV alpha to theta}, \\
	\overline{X}_{i,s}(\bx,t)&=D_{X,s}(\bx) X_i(\bx,t) \label{overX},
\end{align}
and the resulting $\overline{K}_{j,s}(\bk,t)$ are listed in \cref{tab:alpha_methods}.

At the continuous level, two approaches yield identical results and no truncation is required. 
At the discrete level, two approaches require truncation of the infinite integral and exhibit almost no difference in results. 
However, their applicable situations differ. 
The $\mathcal{Y}$-truncation requires the low-rank components $K_j$ to have a Fourier transform. 
This is typically guaranteed by the rapid decay of $K_j$ in the $\bk$-space, which is in turn induced by the rapid decay of the function $f$ in the $\bk$-space. 
In contrast, the $\mathcal{K}$-truncation requires the function $D_{Y,s}(\by)$ to have a Fourier transform, which is not always guaranteed (for example, $D_{Y,1}(\by)=\by$ in the harmonic oscillator).

\begin{table}[h!]
	\centering
	\caption{The resulting $\overline{K}_{j,s}(\bk,t)$ corresponding to truncation methods.}
	\renewcommand{\arraystretch}{2} 
	\begin{tabular}{|c|c|}
		\hline
		Truncation & $\overline{K}_{j,s}(\bk,t)$\\
		\hline
		  $\mathcal{Y}$  &$\left\{\begin{aligned}
			& Y_j(\by,t)=\frac{1}{(2 \pi)^d} \int_{\mathbb{R}^d}  \mathrm{e}^{\mathrm{i} \bk^{\prime} \cdot \by} K_j\left(\bk^{\prime}, t\right) \mathrm{d} \bk^{\prime}, \\
			& \overline{K}_{j,s}(\bk,t)=\frac{1}{\mathrm{i} \hbar}\int_{\mathbb{R}^d} \mathrm{e}^{-\mathrm{i} \bk\cdot\by} Y_j(\by,t) D_{Y,s}(\by) \mathrm{d} \by.
		\end{aligned} \right.$  \\
		\hline
		 $\mathcal{K}$ & $\left\{\begin{aligned}
			& D_{K,s}(\bk)=\frac{1}{(2\pi)^d}\int_{\mathbb{R}^d}\me^{-\mi \bk \cdot \by} D_{Y,s}(\by)\D \by, \\
			& \overline{K}_{j,s}(\bk,t)=\frac{1}{\mi\hbar}\int_{\mathbb{R}^d}K_j(\bk^{\prime},t)D_{K,s}(\bk-\bk^{\prime}) \D \bk^{\prime}.
		\end{aligned} \right.$  \\
		\hline
	\end{tabular}
	\label{tab:alpha_methods}
\end{table}

\section{Dynamical low-rank approximation for the Wigner equation}
\label{DLRA for wigner}
In this section, we substitute \cref{DV alpha to theta} into three substep equations \cref{Z-step}-\cref{L-step} to derive more precise forms of each part, and design their numerical scheme.

\subsection*{$Z$-step} Plugging \cref{DV alpha to theta} into \cref{Z-step} to obtain
\begin{align}
	\frac{\partial}{\partial t}Z_j(\bx,t)= \sum_{s=1}^{R}\sum_{j'=1}^rD_{X,s}(\bx)c_{j,j',s}Z_{j'}(\bx,t)-\frac{\hbar}{m}\sum_{j'=1}^r \ba_{j,j'}\cdot \nabla_{\bx} Z_{j'}(\bx,t). \label{Z}
\end{align}
Here we define time-independent coefficients
\begin{align}
	&A_u\in \mathbb{R}^{r\times r},\quad \quad (A_u)_{i,j}=(\ba_{i,j})_u=\left\langle k_u K_{i},K_j \right\rangle_{\bk}, \quad (u=1,\dots,d), \label{coeff-a}\\
	&C_{s}\in \mathbb{R}^{r\times r},\quad \quad (C_{s})_{i,j}=c_{i,j,s}=\left\langle K_{i},\overline{K}_{j,s}\right\rangle_{\bk}, \quad(s=1,\dots,R). \label{coeff-c}
\end{align}
In this $Z$-step, the coefficients $\ba_{j,j'}$ and $c_{j,j',s}$ are computed first, and equation \cref{Z} is then solved until the next time step.  
Subsequently, a standard QR decomposition is applied to $Z_j$ to obtain new orthonormal  $X_i$ and updated $S_{ij}$ in the end. 
This completes the first $Z$-step.
\subsection*{$S$-step} Plugging \cref{DV alpha to theta} into \cref{S-step} to obtain
\begin{align}
	\frac{\partial}{\partial t} S_{ij}(t)=\frac{\hbar}{m}\sum_{i',j'=1}^r \ba_{j,j'} \cdot \bb_{i,i'}S_{i'j'}(t)-\sum_{s=1}^R\sum_{i',j'=1}^r c_{j,j',s} d_{i,i',s} S_{i'j'}(t). \label{S}
\end{align}
Here we can use time-independent coefficients $\ba_{j,j'}$ and $c_{j,j',s}$ defined in \cref{coeff-a} and \cref{coeff-c} combined with 
\begin{align}
	&B_u\in \mathbb{R}^{r\times r},\quad \quad (B_u)_{i,j}=(\bb_{i,j})_u=\left\langle X_i, \nabla_{x_u} X_j \right\rangle_{\bx}, \quad (u=1,\dots,d), \label{coeff-b}\\
	&D_{s}\in \mathbb{R}^{r\times r},\quad \quad (D_{s})_{i,j}=d_{i,j,s}=\left\langle X_{i},\overline{X}_{j,s} \right\rangle_{\bx},\quad(s=1,\dots,R). \label{coeff-d}
\end{align}
In the $S$-step, since the basis functions $\{K_j(\bk, t) \}_{1\le j \le r}$ remain unchanged during both the $Z$-step and $S$-step, the precomputed values of $\ba_{j,j'}$ and $c_{j,j',s}$ from the previous $Z$-step can be used directly.
The newly introduced coefficients $\bb_{i,i'}$ and $d_{i,i',s}$ should be computed first, the $S$-step equation \cref{S} can then be solved to obtain the updated $S$.
This completes the second $S$-step.

\subsection*{$L$-step} Plugging \cref{DV alpha to theta} into \cref{L-step} to obtain
\begin{align}
	\frac{\partial}{\partial t} L_i(\bk, t)=\sum_{s=1}^R\sum_{i',j'=1}^r d_{i,i',s} S_{i'j'}(t)\overline{K}_{j',s}(\bk,t)  - \frac{\hbar}{m}\sum_{i'=1}^r \bk \cdot \bb_{i,i'} L_{i'}(\bk,t), \label{L}
\end{align}
with coefficients $\bb_{i,i'}$ and $d_{i,i',s}$ have defined in \cref{coeff-b} and \cref{coeff-d}. 
Since the basis functions $\{X_i(\bx, t) \}_{1\le i \le r}$ remain unchanged during $S$-step, the precomputed values of $\bb_{i,i'}$ and $d_{i,i',s}$ from the previous $S$-step can be used directly to solve \cref{L}.
Then apply standard QR decomposition to $L_i$ to get new orthonormal  $K_j$ and updated $S_{ij}$ in the end. 
This completes the last $L$-step.

\begin{remark}
	Coefficient matrices $A_u$ and $D_{s}$ are symmetric by construction, which ensures hyperbolicity of whole system.
\end{remark}

A detailed account of first order accurate DLRA for the Wigner equation is given in \cref{alg:brief first DLRA for wigner}.
The solutions corresponding to the $Z$-step, the $S$-step and the $L$-step are combined by a first order Lie splitting scheme in order to recover $f(\bx, \bk, t)$. 
The sequence of these three steps is not unique \cite{lubich2014projectorsplitting}, but we only consider the sequence of $ZSL$ in this work.

\begin{algorithm}
	\caption{A first-order accurate dynamical low-rank algorithm for the Wigner equation}
	\renewcommand{\algorithmicrequire}{\textbf{Input:}}
	\renewcommand{\algorithmicensure}{\textbf{Output:}}
	\label{alg:brief first DLRA for wigner}
	\begin{algorithmic}[1]
		\REQUIRE $D_{X,s},D_{Y,s},X_i^0,S_{ij}^0,K_j^0$ such that $f(\bx, \bk, 0)= \sum_{i,j=1}^{r} X_i^0(\bx)S_{ij}^0 K_j^0(\bk)$ with time step $\tau$.
		\ENSURE $X_i^1,S_{ij}^3,K_j^1$ such that $f(\bx, \bk, \tau)= \sum_{i,j=1}^{r} X_i^1(\bx)S_{ij}^3 K_j^1(\bk)$.
		\STATE Compute $Z_i^0$ using $X_i^0$ and $S_{ij}^0$.
		\STATE Compute $\overline{K}_{j,s}$ using $D_{X,s},D_{Y,s},$ and $K_j^0$ according to \cref{tab:alpha_methods}.
		\STATE Compute $\ba_{j,j'}$ and $c_{j,j',s}$ using $K_j^0$ and $\overline{K}_{j,s}$.
		\STATE Solve \cref{Z} with initial value $Z_i^0$ until time $\tau$ to obtain $Z_i^1$.
		\STATE Perform a QR decomposition of $Z_i^1$ to obtain $X_i^1$ and $S_{ij}^1$.
		\STATE Compute $\bb_{i,i'}$ and $d_{i,i',s}$ using $X_j^0$.
		\STATE Solve \cref{S} with initial value $S_{ij}^0$ until time $\tau$ to obtain $S_{ij}^2$.
		\STATE Compute $L_i^0$ using $K_i^0$ and $S_{ij}^2$.
		\STATE Solve \cref{L} with initial value $L_i^0$ until time $\tau$ to obtain $L_i^1$.
		\STATE Perform a QR decomposition of $L_i^1$ to obtain $K_i^1$ and $S_{ij}^3$.
	\end{algorithmic}
\end{algorithm}

We have restricted our analysis to first-order Lie splitting. In principle, extending the method to second-order Strang splitting is straightforward; further implementation details can be found in \cite{einkemmer2018lowrank,cassini2022efficient}.
Although both first-order Lie splitting and second-order Strang splitting have been observed to achieve their respective convergence orders in practice, existing theoretical analyses can only guarantee first-order accuracy \cite{kieri2016discretized}. This limitation arises because the projection operator lacks Lipschitz continuity \cite{koch2007dynamical}.
To date, this phenomenon remains not fully understood theoretically \cite{ceruti2024robust}.
The numerical results presented in this work indicate that, even when employing a first-order Lie splitting method, second-order accuracy can still be attained provided that each splitting substep itself achieves second-order accuracy.
This numerical observation warrants further theoretical investigation in future work. 

\subsection{Discretization}
\label{sec:Discretization}
\label{subsec:Discretization}
For the purpose of numerical solution, we first truncate the computational domain to a bounded rectangular region $\Omega=\Omega_{\bx}\times \Omega_{\bk},$ then we consider the uniform grid mesh in $\bx$-space and $\bk$-space, $\bx$-space using $N_x$ grid points in every coordinate direction with the spacing $\Delta x_u$ and $\bk$-space using $N_k$ grid points in every coordinate direction with the spacing $\Delta k_u.$ 
For the dual $\by$-space, we choose $N_y$ grid points in every coordinate direction with the spacing $\Delta y_u$. 

Once the grid discretization is determined, for a fixed time $t$, we can obtain $X(t)\in \mathbb{R}^{N_{x}^d\times r}$ to be the matrix composed of row vectors collected from $X_i$ on the chosen spatial grid. 
It's obvious that the resulting matrix entries depend on both the discretization performed and the ordering of the grid points. 
Similarly, $K(t)\in \mathbb{R}^{N_{k}^d\times r}$ can be defined in the same way. 
Then, the low-rank approximation of the Wigner function can be reconstructed at any time from the low-rank factors via the computation
\begin{align}
	F(t)=X(t)S(t)K^T(t) \in \mathbb{R}^{N_{x}^d\times N_{k}^d}.
\end{align}
However, due to its prohibitive cost, this reconstruction is unnecessary for our algorithm.
We only employ it in combination with Singular Value Decomposition (SVD) \cite{eckart1936approximation} to obtain the initial values of the low-rank factors and to reconstruct the solution at selected critical time points. 
Similarly, the separable assumption \cref{DV low-rank} for the difference potential $D_V$ can be computed via SVD, and this computation is required only once during the initialization stage. 
In each time step, since a uniform grid is employed, $\overline{K}_{j,s}$ can be efficiently computed via the standard Fast Fourier transform (FFT). 
It is also worth noting that the Poisson summation formula is often used in conjunction with truncation \cite{chen2022higherorder,chen2019high,chen2019numerical,xiong2016advectivespectralmixed}, however, it is not required in the implementation of our DLRA algorithm. 

\subsection{Numerical scheme}
In this subsection, we provide a detailed description of the procedure for solving \cref{Z}, \cref{S}, and \cref{L}. 

We first consider the $S$-step \cref{S}. It is an ordinary differential equation and there is no source of stiffness in it. 
Hence, we perform its time integration by classical explicit fourth order Runge-Kutta scheme RK4.

Next, let us consider the $Z$-step \cref{Z} and use the vector notation $Z=(Z_1,...,Z_r)$, we can rewrite the evolution equation for the $Z$-step \cref{Z} in matrix formulation as follows:
\begin{align}
\frac{\partial}{\partial t}Z(\bx,t)&=\sum_{s=1}^R D_{X,s}(\bx) Z(\bx,t)C_s^T
-\frac{\hbar}{m}\sum_{u=1}^d \nabla_{x_u}Z(\bx,t) A_u^T . \label{MA-Z}
\end{align}
For simplicity, we only consider the two-dimensional case here. 
By further splitting the $Z$-step, we obtain the following two subflows:
\begin{align}
	\frac{\partial}{\partial t}Z(\bx,t)&= -\frac{\hbar}{m} \nabla_{x_1}Z(\bx,t) A_1^T , \label{MA_Y_1}\\
	\frac{\partial}{\partial t}Z(\bx,t)&= -\frac{\hbar}{m} \nabla_{x_2}Z(\bx,t) A_2^T+\sum_{s=1}^R D_{X,s}(\bx) Z(\bx,t)C_s^T. \label{MA_Y_2}
\end{align}
As we mentioned before, the coefficients matrix $A_1$ and $A_2$ are symmetric, thus equations \cref{MA_Y_1} and \cref{MA_Y_2} are hyperbolic. 
It is evident that \cref{MA_Y_1} can be solved in time by the method of characteristics, while \cref{MA_Y_2} can be efficiently solved via an exponential integrator (see \cite{hochbruck2010exponential} for more details).

Finally, we consider the $L$-step \cref{L}. 
Extending the vector notation similar to $Z$ to $X, K, Y, L, S,$ and $\overline{K}_s$, we can obtain 
\begin{align}
	\frac{\partial}{\partial t}L(\bk,t)& =\sum_{s=1}^R \overline{K}_s S^T D_{s}^T
	- \frac{\hbar}{m} \sum_{u=1}^{d} k_u  L(\bk,t) B_u^T. \label{MA-L}
\end{align}
Unlike the $Z$-step, the hyperbolicity is maintained by the coefficient matrix $D_s$, hence it is necessary to examine the relationship between $L$ and $\overline{K}_s S^T$ and to exploit this connection in solving the $L$-step equation. 

For example, for $\mathcal{Y}$-truncation, using the relationship provided in \cref{tab:alpha_methods}, apply the Fourier transform in $\bk$-direction to \cref{MA-L} and define
\begin{align}
	L(\bk,t)=K(\bk,t)S^T&=\int_{\mathbb{R}^d} \mathrm{e}^{-\mathrm{i} \bk\cdot\by} U(\by,t)\mathrm{d} \by, \label{L to U}\\
	\frac{\hbar}{m} \sum_{s=1}^{d} k_s  L(\bk,t) B_s^T&=\int_{\mathbb{R}^d} \mathrm{e}^{-\mathrm{i} \bk\cdot\by} V(\by,t) \mathrm{d} \by, \label{L to V}
\end{align}
then the $L$-step \cref{MA-L} turns to
\begin{align}
	\frac{\partial}{\partial t}U(\by,t)=\frac{1}{\mathrm{i} \hbar}\sum_{s=1}^{R} D_{Y,s}(\by)U(\by,t)D_{s}^T-V(\by,t). \label{MA-L-Y}
\end{align}
Therefore, we proceed to solve equation \cref{MA-L-Y} involving the variable $U(\by,t)$, where $V(\by,t)$ can be regarded as a function of $U(\by,t)$.
When $R$ is approximately the same size as the dimension $d$, similar to the situation in the $Z$-step, we can likewise exploit the symmetry of the coefficient matrix $D_{s}$ and combine a splitting method with an exponential integrator to solve the corresponding equations.
Conversely, when $R$ is relatively large, the Runge-Kutta method is more suitable for \cref{MA-L-Y}. Similarly, for $\mathcal{K}$-truncation, utilizing the convolution property of Fourier transform and the notation in \cref{L to U} and \cref{L to V}, we can obtain the evolution equation for $U(\by,t)$ as same as \cref{MA-L-Y}. 

\subsection{Complexity analysis} 
Then we considered the time complexity of the proposed \cref{alg:brief first DLRA for wigner}. Here, we consider the case of $N_y=\mathcal{O}(N_k).$ 

We first consider the computational complexity required for calculating $\overline{K}_{j,s}(\bk,t)$. 
The additional storage requirements are all of order $\mathcal{O}(N_k^d(R+r))$ and the required time complexity is $\mathcal{O}(N_k^d\log(N_k)Rr)$. 
Then consider the other steps in \cref{alg:brief first DLRA for wigner}, the storage cost and time complexity for computing the coefficients and solving the equations is shown in \cref{tab:alg_complexity}. 
For completeness, the time complexity of both QR decomposition and matrix multiplication can be summarized as $\mathcal{O}( (N_k^d+N_x^d) r^2 ).$
In summary, compared to the full-grid method, which exhibits a time and space complexity of $\mathcal{O}(N_x^d N_k^d),$ our DLRA method offers significant advantages.

\begin{table}[h]
	\centering
	\renewcommand\arraystretch{1.4}
	\caption{Complexity analysis of the main steps in the \cref{alg:brief first DLRA for wigner}.}
	\label{tab:alg_complexity}
	\begin{tabular}{| c | c | c |}
		\hline
		Steps  & Storage cost & Time complexity\\ \hline
		$A_u$ & $\mathcal{O}(r^2)$ & $\mathcal{O}(N_k^d r^2)$  \\ \hline
		$B_u$ &  $\mathcal{O}(r^2)$ & $\mathcal{O}(N_x^d \log(N_x)r+N_x^d r^2)$ \\ \hline
		$C_s,D_s$ & $\mathcal{O}(R r^2)$ & $\mathcal{O}(N_k^d R r^2)$ \\ \hline
		$Z$-step & $\mathcal{O}(N_x^d r)$ &  $\mathcal{O}(N_x^d \log(N_x) r+N_x^d R r^2)$  \\ \hline
		$S$-step & --- &  $\mathcal{O}(Rr^4)$  \\ \hline
		$L$-step & $\mathcal{O}(N_k^d r)$ &  $\mathcal{O}(N_k^d \log(N_k) r+ N_k^d R r^2)$  \\ \hline
	\end{tabular}
\end{table}

\begin{remark}
	Due to the convolution involved, $\mathcal{K}$-truncation requires zero-padding of factor $2$ to avoid aliasing errors; however, its complexity does not differ from that of the $\mathcal{Y}$-truncation.
\end{remark}

\section{Numerical experiments}
\label{Numerical experiments}
In this section, we evaluate the accuracy and efficiency of the proposed DLRA method. 
As we mentioned in \Cref{Separated approximation for pdo}, we adopt the $\mathcal{Y}$-truncation in the subsequent numerical experiments.

We assess the accuracy of our algorithm using a series of typical benchmark tests for the Wigner equation, including the quantum harmonic oscillator and Gaussian barrier scattering. These test cases are characterized by non-stationary and non-low-rank features, which are significant for extending the applicability of the DLRA method. 
Once numerical accuracy has been verified, we can investigate typical quantum systems such as the electron-electron scattering and Helium-like system, and reveal fundamental quantum effects such as the Pauli exclusion principle and the Heisenberg uncertainty principle.

The numerical performance is evaluated using the  $L^2$ error $\varepsilon_2(t)$ and the  $L^{\infty}$ error $\varepsilon_{\infty}(t)$, defined as follows: 
\begin{align}
	\varepsilon_{2}(t)&= \left[\iint_{\mathcal{X}\times \mathcal{K}} \left(f^{\textup{ref}}\left(\bm{x},\bm{k},t\right)-f^{\textup{num}}\left(\bm{x},\bm{k},t\right)\right)^{2}\textup{d}\bm{x}\textup{d} \bm{k}\right]^{\frac{1}{2}}, \\
	\varepsilon_{\infty}(t) &=\max_{(\bm{x},\bm{k})\in\mathcal{X}\times \mathcal{K}}\big |f^{\textup{ref}}\left(\bm{x},\bm{k},t\right)-f^{\textup{num}}\left(\bm{x},\bm{k},t\right) \big |,
\end{align}
where $f^{\textup{ref}}$ and $f^{\textup{num}}$ denote the reference and DLRA numerical solution, respectively.

For 4D problem, we adopt the reduced Wigner function onto the 2D phase-space plane $(x_j,k_j)$,
\begin{align}
	W_s(x,k,t)=\iint_{\mathbb{R}\times\mathbb{R}} f(\bm{x},\bm{k},t) \D \bx_{3-s}\D \bk_{3-s}, \quad s=1,2,
\end{align}
or reduced one-body Wigner function \cite{cancellieri2007effect}
\begin{align}
	W(x,k,t)=W_1(x,k,t)+W_2(x,k,t).
\end{align}
In practice, the low-rank assumption of $f(\bx,\bk,t)$ in \cref{low-rank} yields
\begin{align}
	W_s(x,k,t)=\sum_{i,j=1}^rS_{ij}(t)\int_{\mathbb{R}} X_i(\bx,t) \D \bx_{3-s} \int_{\mathbb{R}} K_j(\bk,t) \D \bk_{3-s},
\end{align}
so that $W_s$ can be obtained directly from the low‑rank components $X_i,S_{ij}$ and $K_j$ without first recovering the full-grid numerical solution $f^{\textup{num}}$.

\subsection{Quantum harmonic oscillator}
\label{subsec: harmonic}
In this subsection, we consider 2D quantum harmonic oscillator $V(x)=m\omega x^2/2$ with $R=1$.
The exact solution can be obtained as $f(x,k,t)=f_0(x(t),k(t))$, where $(x(t),k(t))$ are the characteristic lines given by
\begin{equation}
	\left\{\begin{aligned}
		& x(t)=\cos (\sqrt{\omega} t) x(0)-\frac{\hbar}{m \sqrt{\omega}} \sin (\sqrt{\omega} t) k(0), \\
		& k(t)=\frac{m \sqrt{\omega}}{\hbar} \sin (\sqrt{\omega} t) x(0)+\cos (\sqrt{\omega} t) k(0) .
	\end{aligned}\right.
\end{equation}
For simplicity, we still choose parameters $\hbar=m=\omega=1$, equip with an initial standard Gaussian wave packet (GWP)
\begin{align}
	f_0(x,k)=\frac{1}{\pi}\me^{-\frac{(x-x_0)^2}{2\sigma^2}-2\sigma^2(k-k_0)^2}, \label{GWP}
\end{align}
where $x_0=1$ is the center of GWP and $k_0=0$ is the initial wave number. 
Then reference solution can be written as
\begin{align}
	f(x,k,t)=f_0(\cos(t)x-\sin(t)k, \sin(t)x+\cos(t)k), 
\end{align}
which means that reference solution has a period of $2\pi$, thus we choose final time $T=2\pi$. 
The computational domain is $\Omega=\Omega_x \times \Omega_k=[-12,12]\times[-6.4,6.4]$ with $N_x=N_k=128$. 
Since we mainly focus on the convergence with respect to time step $\Delta t = T/N_t$ and rank $r$, several groups of simulations under $N_t=5000, 10000, 20000, 40000, 80000$ and $r=5,7,9,10,15,20,25,30$ are performed.

\begin{figure}[h!]
	\centering
	\subfigure[$\|\cdot\|_2$ for best rank-r approximation. \label{fig::harmonic_2norm}]
	{\includegraphics[width=0.45\textwidth,height=0.4\textwidth]{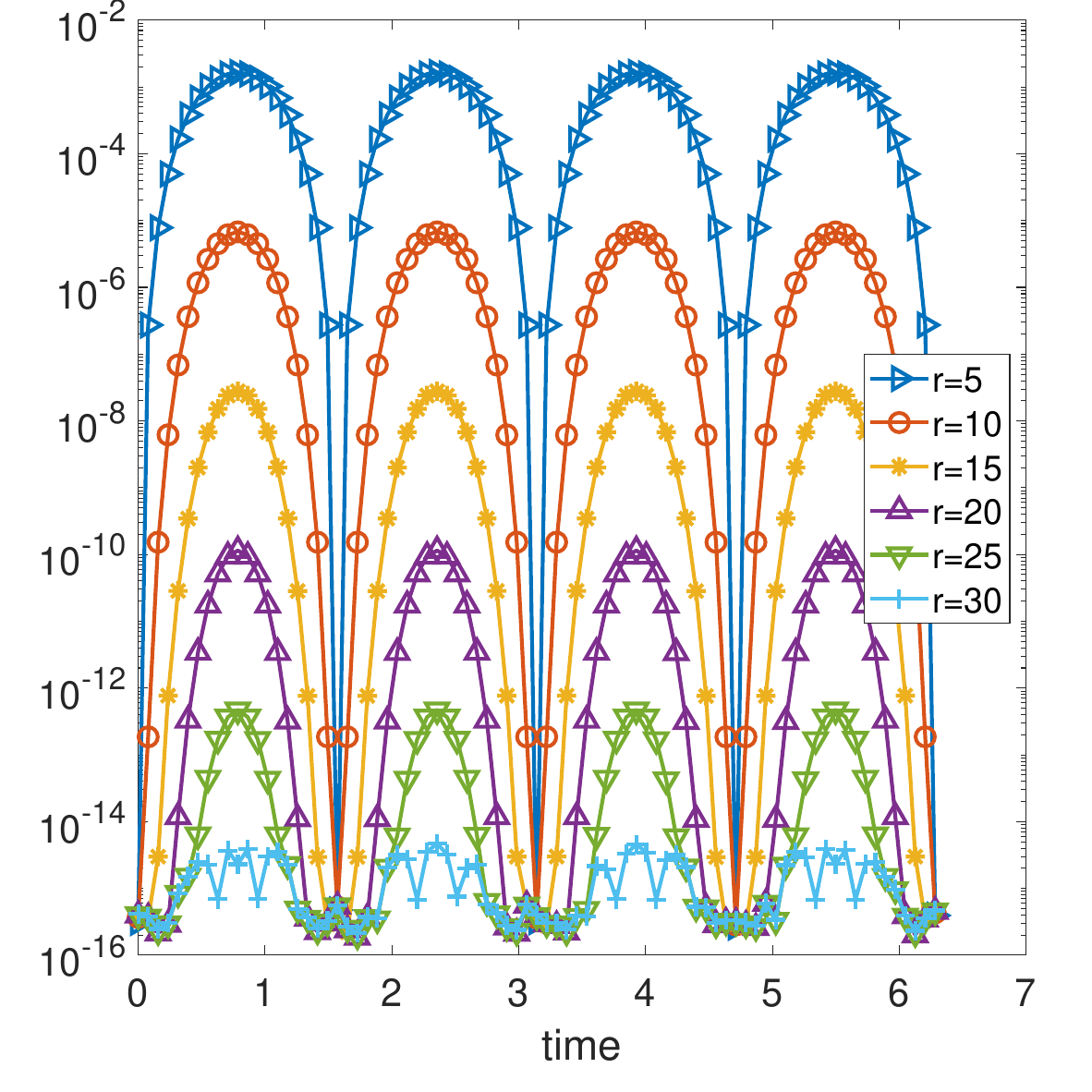}}
	\subfigure[Convergence with respect to $\Delta t$.\label{harmonic_convergence}]
	{\includegraphics[width=0.45\textwidth,height=0.4\textwidth]{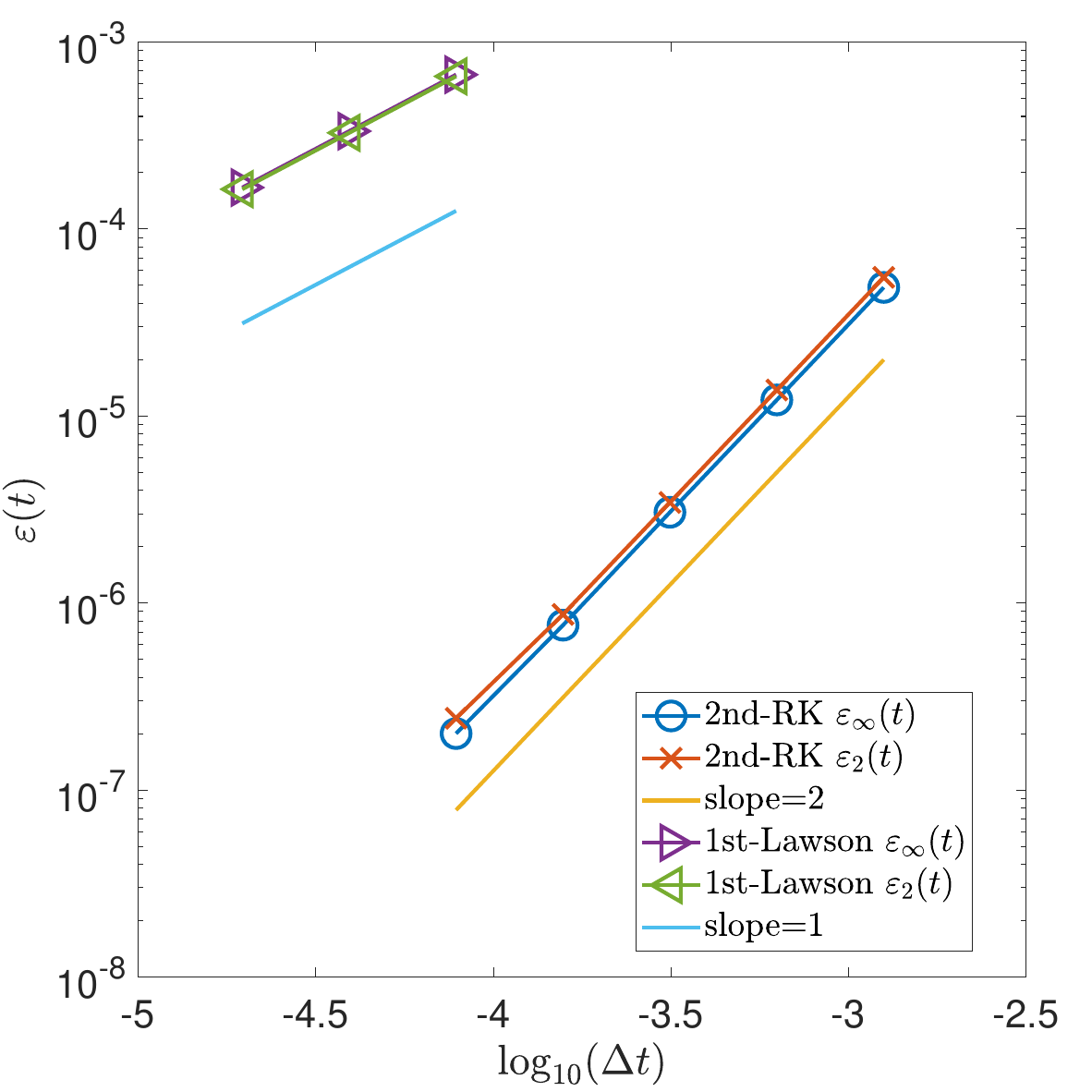}}\\
	\subfigure[$\varepsilon_{2}(t)$ under $N_t=10000$ and different $r$.\label{fig::harmonic_e2r}]
	{\includegraphics[width=0.45\textwidth,height=0.4\textwidth]{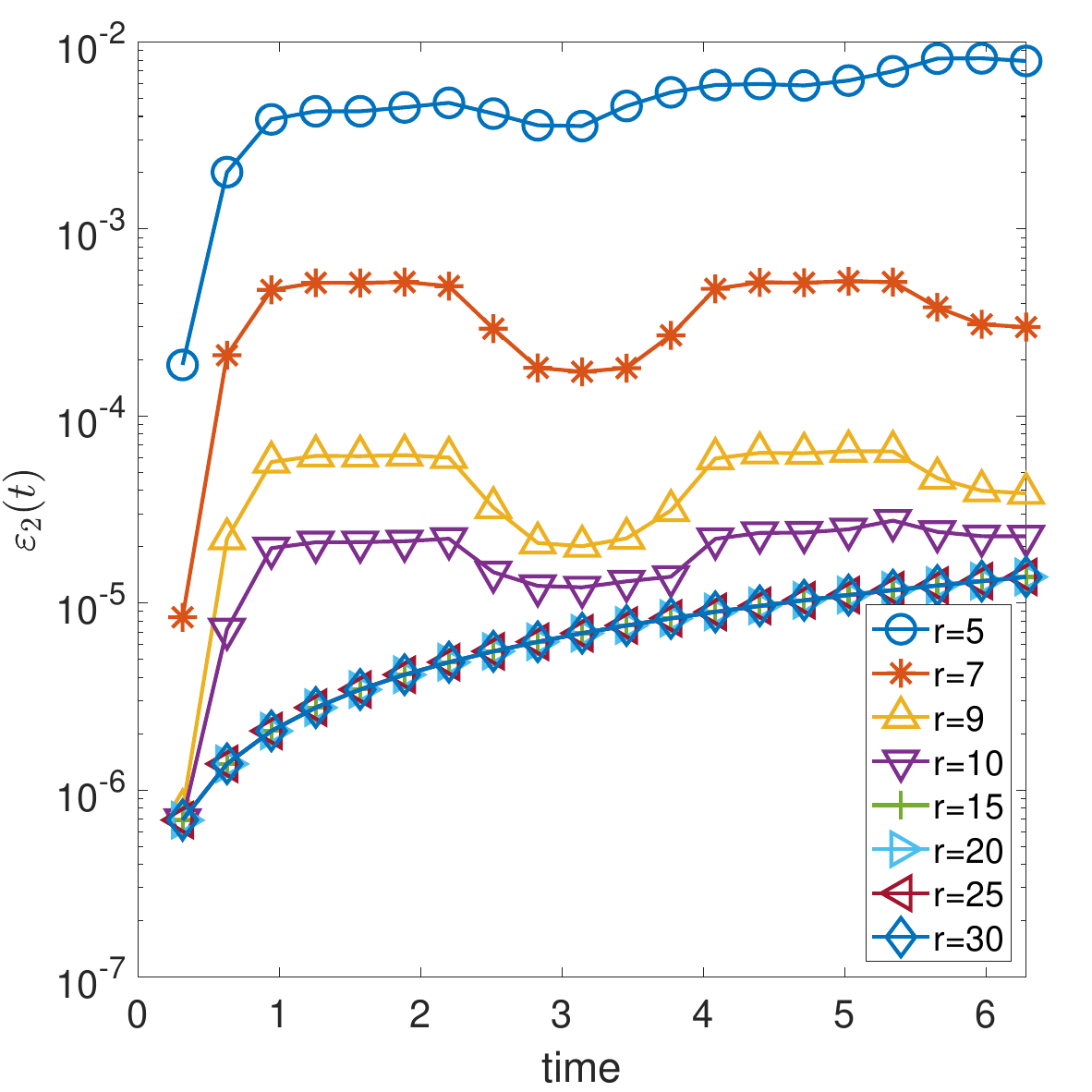}}
	\subfigure[$\varepsilon_{\infty}(t)$ under $N_t=10000$ and different $r$.\label{fig::harmonic_eir}]
	{\includegraphics[width=0.45\textwidth,height=0.4\textwidth]{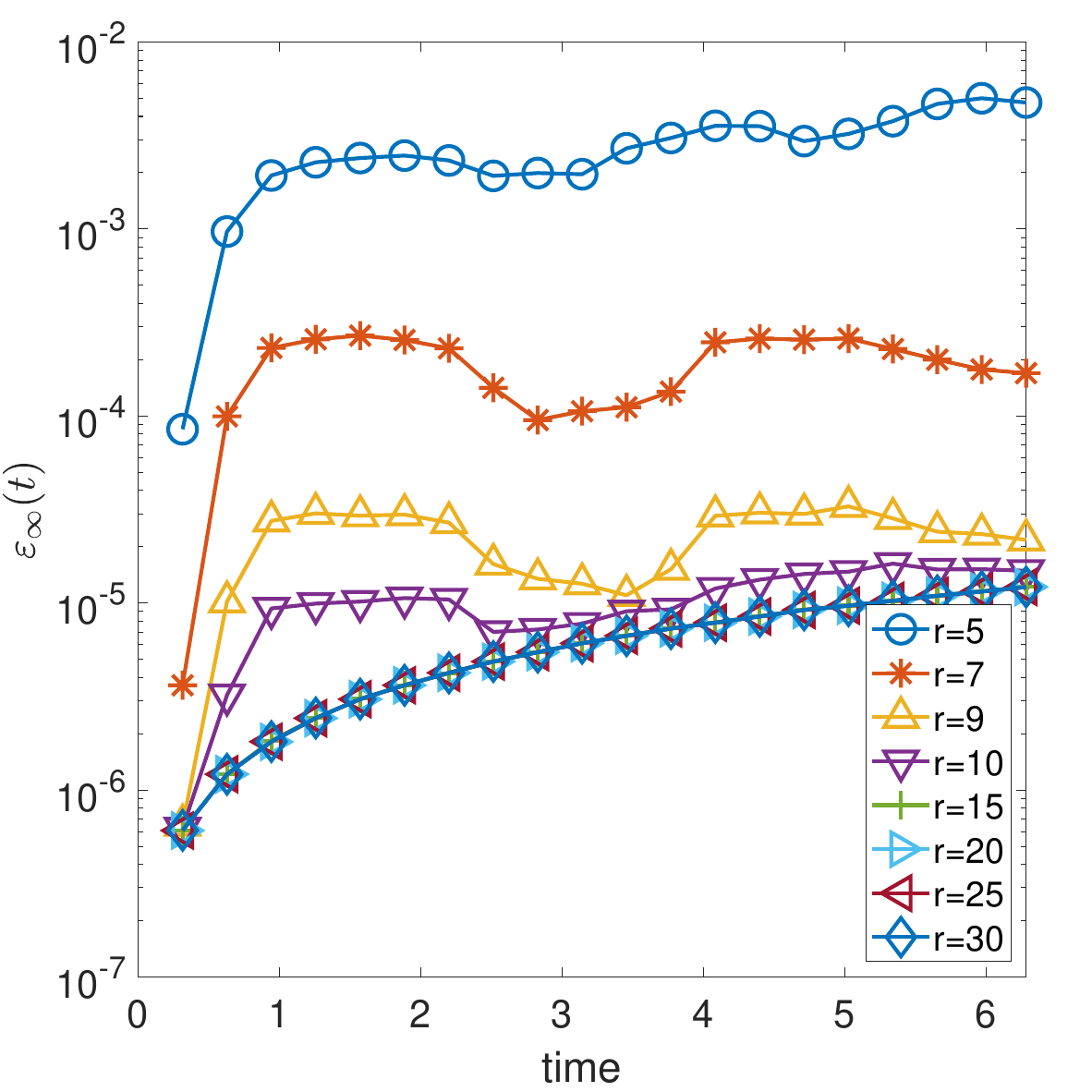}}\\
	\subfigure[$\varepsilon_{2}(t)$ under $r=20$ and different $N_t$.\label{fig::harmonic_e2dt}]
	{\includegraphics[width=0.45\textwidth,height=0.4\textwidth]{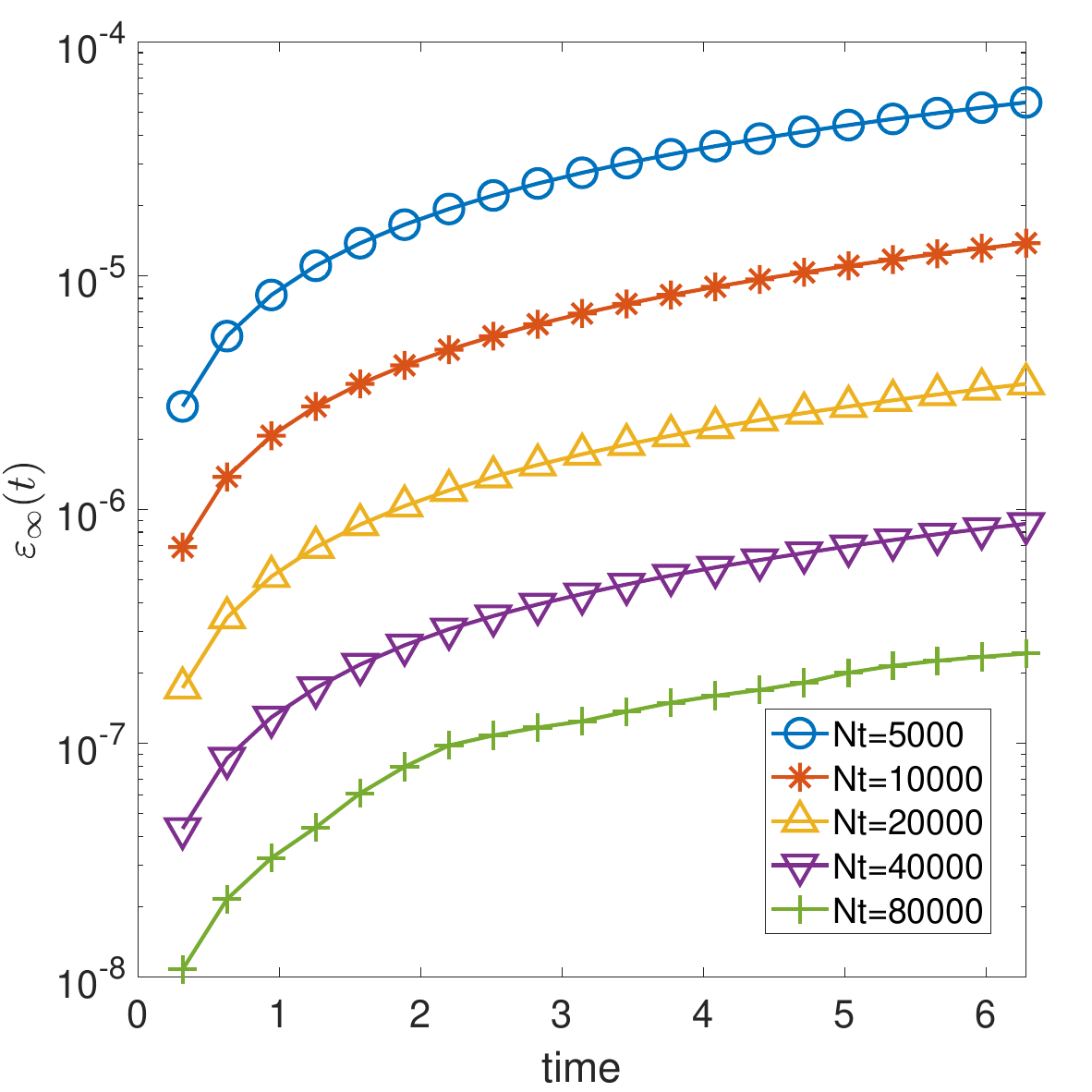}}
	\subfigure[$\varepsilon_{\infty}(t)$ under $r=20$ and different $N_t$.\label{fig::harmonic_eidt}]
	{\includegraphics[width=0.45\textwidth,height=0.4\textwidth]{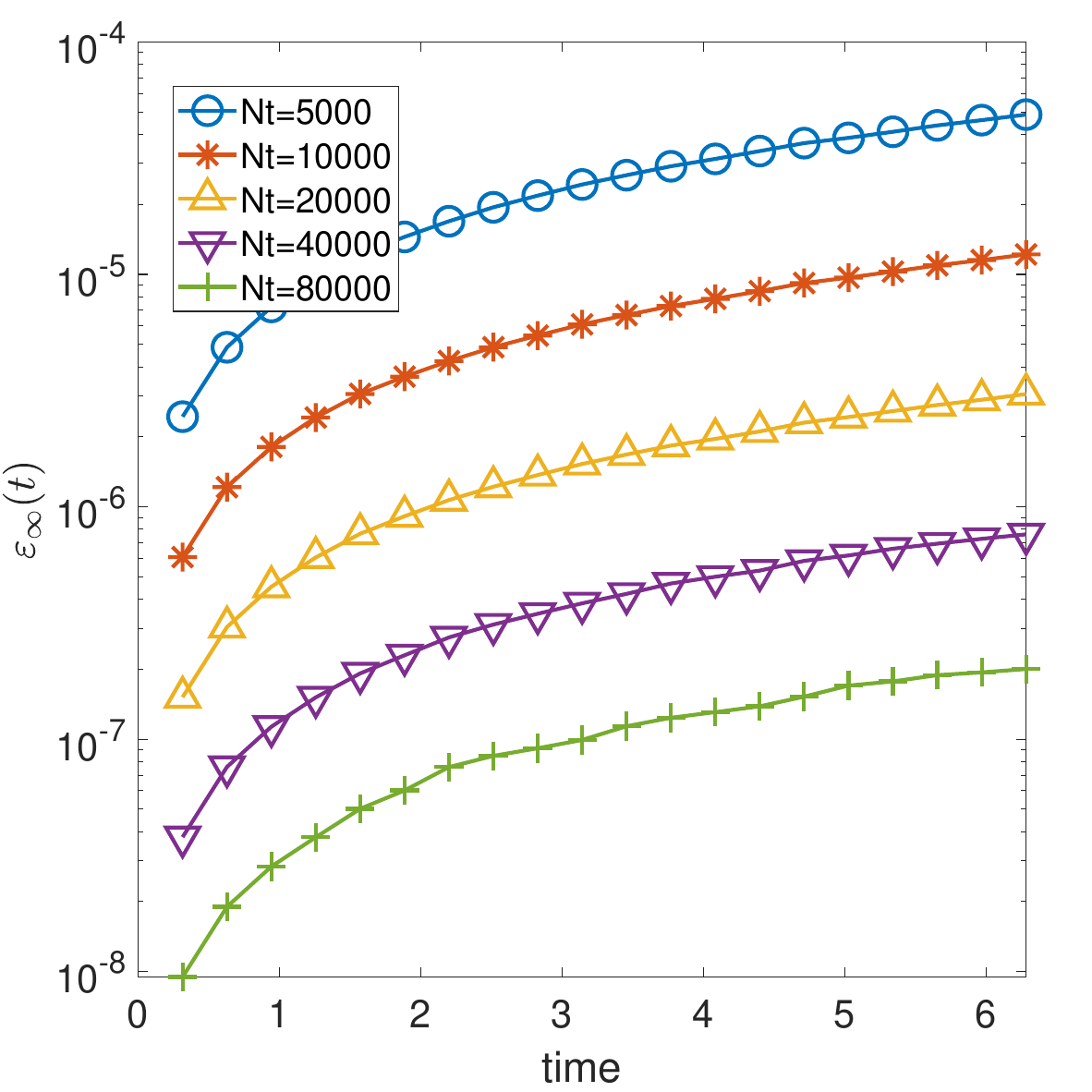}}
	\caption{2D quantum harmonic oscillator.}
	\label{fig::harmonic}
\end{figure}

First, we examine the rank of the reference solution. Once the spatial discretization is fixed, the reference solution can be treated as a matrix, and its effective rank can be estimated via SVD.
The $L_2$ norm of the error between the reference solution matrix and its best rank-$r$ approximation for $r=5,10,15,20,25,30$ is shown in   \cref{fig::harmonic_2norm}, which is equal to $(r+1)$-th singular value of the error matrix.
Therefore, we observe that the singular values of the solution decay in the range of machine accuracy without a distinct gap, except $4$ specific times $t=0,\frac{\pi}{2},\pi,\frac{3\pi}{2}$ in one period, which the exact rank equal to 1. 
Note that the solution from our DLRA algorithm can not be more accurate than the best rank $r$ approximation obtained via SVD. 
Hence, this numerical test either does not exhibit a genuinely low‑rank structure, or its effective rank is relatively high, a behavior that differs from typical DLRA numerical experiments.

Second, we examine the convergence with respect to the rank $r$. 
We fix time discretization parameters $N_t=10000$ and variate the rank $r$, the resulting numerical error $\varepsilon_{2}(t)$ and $\varepsilon_{\infty}(t)$ are shown in the \cref{fig::harmonic_e2r} and \cref{fig::harmonic_eir}, respectively. 
It can be observed that $\varepsilon_{2}(t)$ and $\varepsilon_{\infty}(t)$ behave roughly similarly. 
When $r$ increases from 5 to 15, the numerical error decreases significantly. 
However, as $r$ grows further from 15 to 30, the reduction slows down or even stops, this is because the error at this time is mainly controlled by time step $\Delta t$. 

Last, we investigate the convergence with respect to the time step $\Delta t$. 
Error indicators $\varepsilon_{2}(t)$ and $\varepsilon_{\infty}(t)$ for a fixed rank $r=20$ and different time discretization $N_t$ are shown in the \cref{fig::harmonic_e2dt} and \cref{fig::harmonic_eidt}. 
We observe that the performances of $\varepsilon_{2}(t)$ and $\varepsilon_{\infty}(t)$ are basically the same, they show a clear convergence trend as we increase the time partition parameters $N_t$ and the convergence rate is plotted in \cref{harmonic_convergence}. 
However, \cref{harmonic_convergence} reveals two distinct convergence orders.
The difference stems from the numerical schemes employed for the $Z$-step \cref{MA-Z} and $L$-step \cref{MA-L}:
when these steps are solved with a second‑order Runge–Kutta (2nd‑RK) or a second‑order Lawson (2nd‑Lawson) method, the numerical error shows second‑order convergence;
when a first‑order Lawson (1st‑Lawson) method is used, only first‑order convergence is obtained. 
This indicates that the projector‑splitting error actually exhibits a numerical accuracy no lower than second order, even though a first‑order projector‑splitting algorithm is employed in this setting.

\subsection{Gaussian barrier scattering}
\label{2D Gaussian barrier scattering}
To further validate the overall performance of the DLRA for Wigner equation, we simulate the Gaussian barrier scattering \cite{shao2011adaptive,shao2015comparison,xiong2016advectivespectralmixed} of the GWP to investigate its convergence rate. 
In this numerical experiment, a Gaussian potential barrier is introduced into the spatial domain, given by
\begin{align}
	V(x)=H_P\exp\left(-\frac{(x-x_B)^2}{2w^2}\right),
\end{align}
where $H_P=0.3$ is the height of the potential barrier, $x_B=30$ is the center of the computational domain and $w=1$ is the dispersion of the potential barrier.
We adopt the computational domain $\Omega=\Omega_x \times \Omega_k=[0, 60] \times [-4 ,4]$ with spacing $\Delta y = \pi/4 $.
The initial Wigner function is chosen to be the GWP \cref{GWP} with $x_0=15$, $k_0=0.7$ and $\sigma=2.825$.
And the other parameters are: the reduced Planck constant $\hbar=0.658211899$, the effective mass $m=0.0665m_e$, the
stationary electron mass $m_e=5.68562966$ and the final time $T=20$.
                         
\begin{figure}[h!] 
	\centering
	\subfigure[$\|\cdot\|_2$ for best rank-r approximation. \label{fig::2Dbar_approxr}]
	{\includegraphics[width=0.45\textwidth,height=0.4\textwidth]{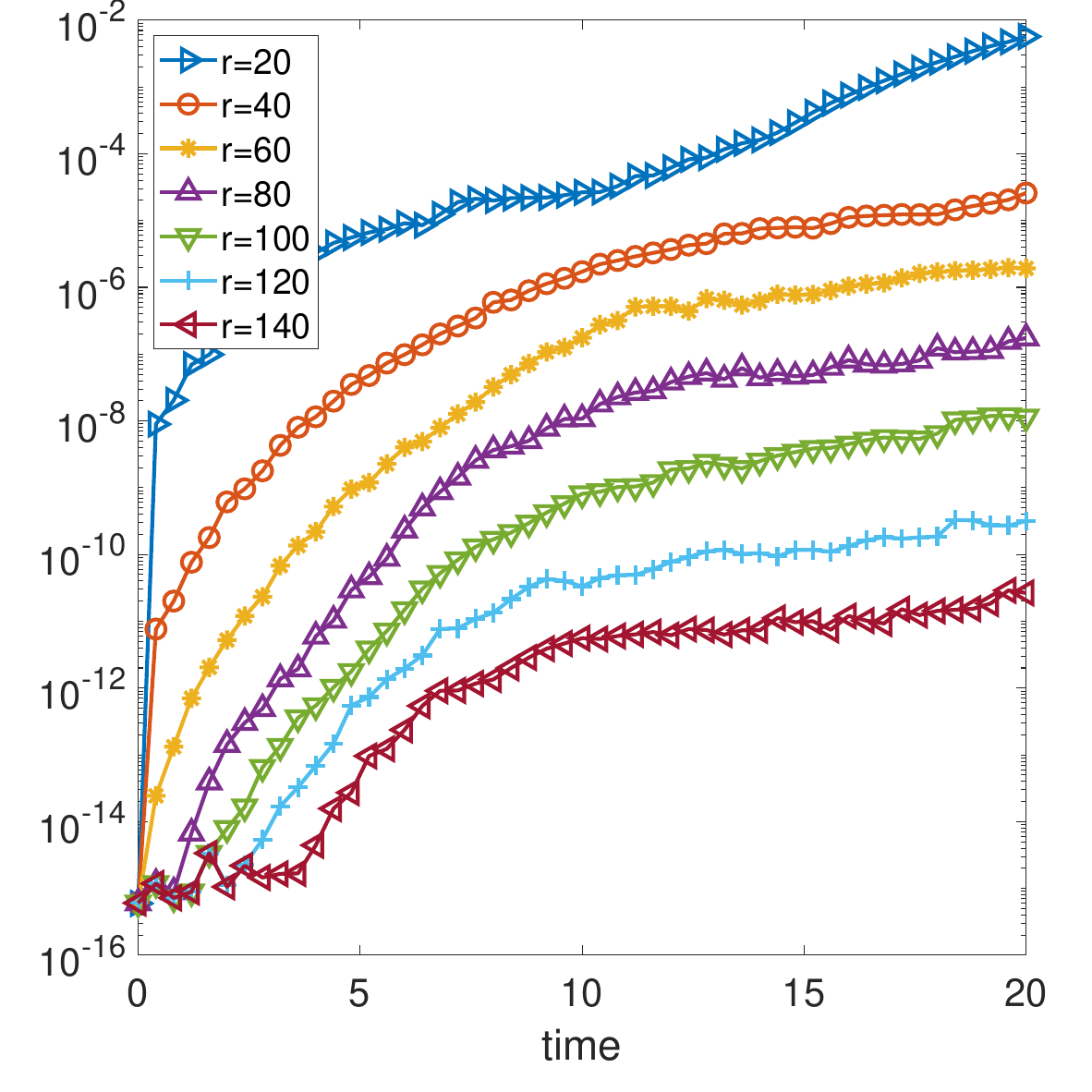}}
	\subfigure[Convergence with respect to $\Delta t$. \label{fig::2Dbar_rate}]
	{\includegraphics[width=0.45\textwidth,height=0.4\textwidth]{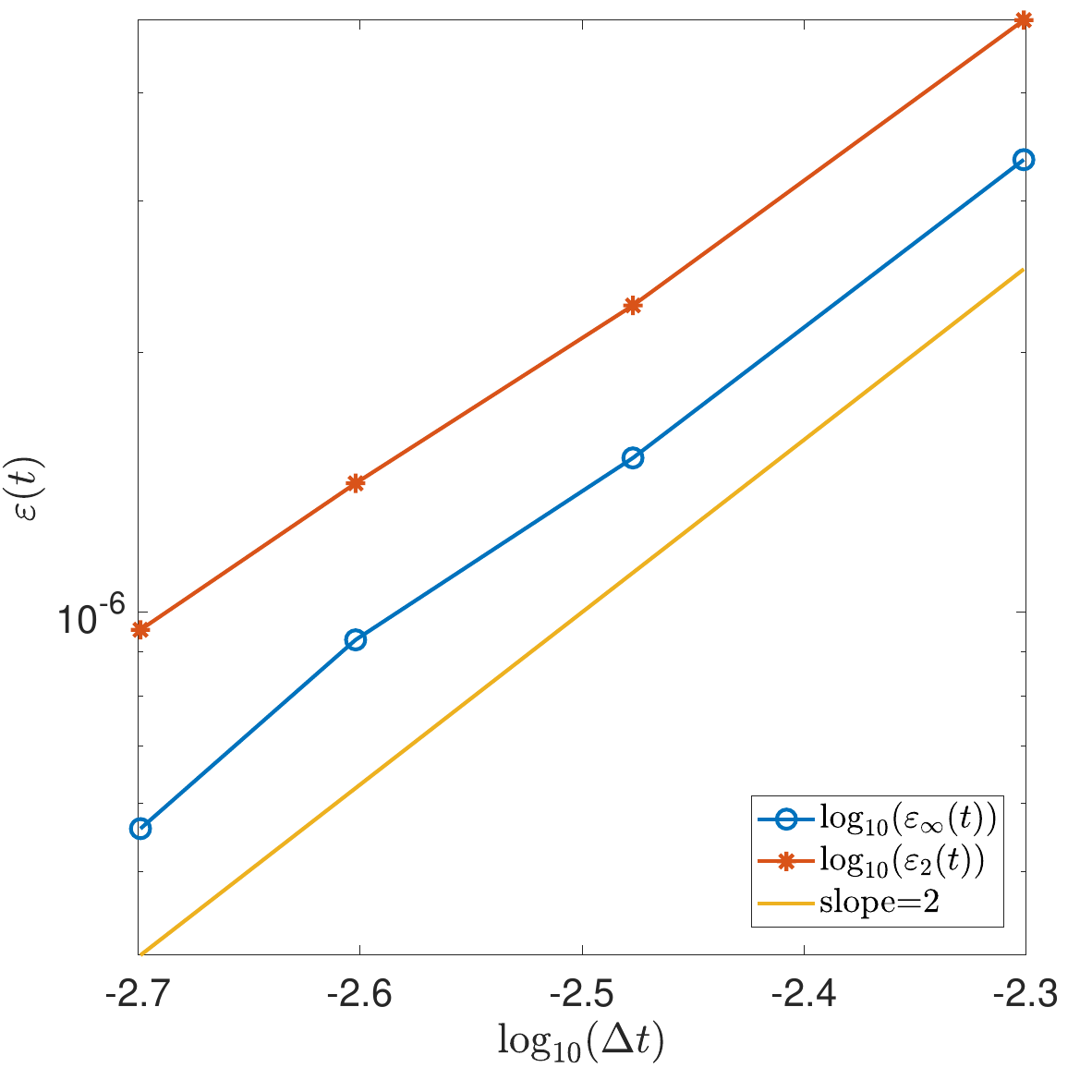}}\\
	\subfigure[$\varepsilon_{2}(t)$ under $N_t=10000$ and different $r$. \label{fig::2Dbar_error_r_e2}]
	{\includegraphics[width=0.45\textwidth,height=0.4\textwidth]{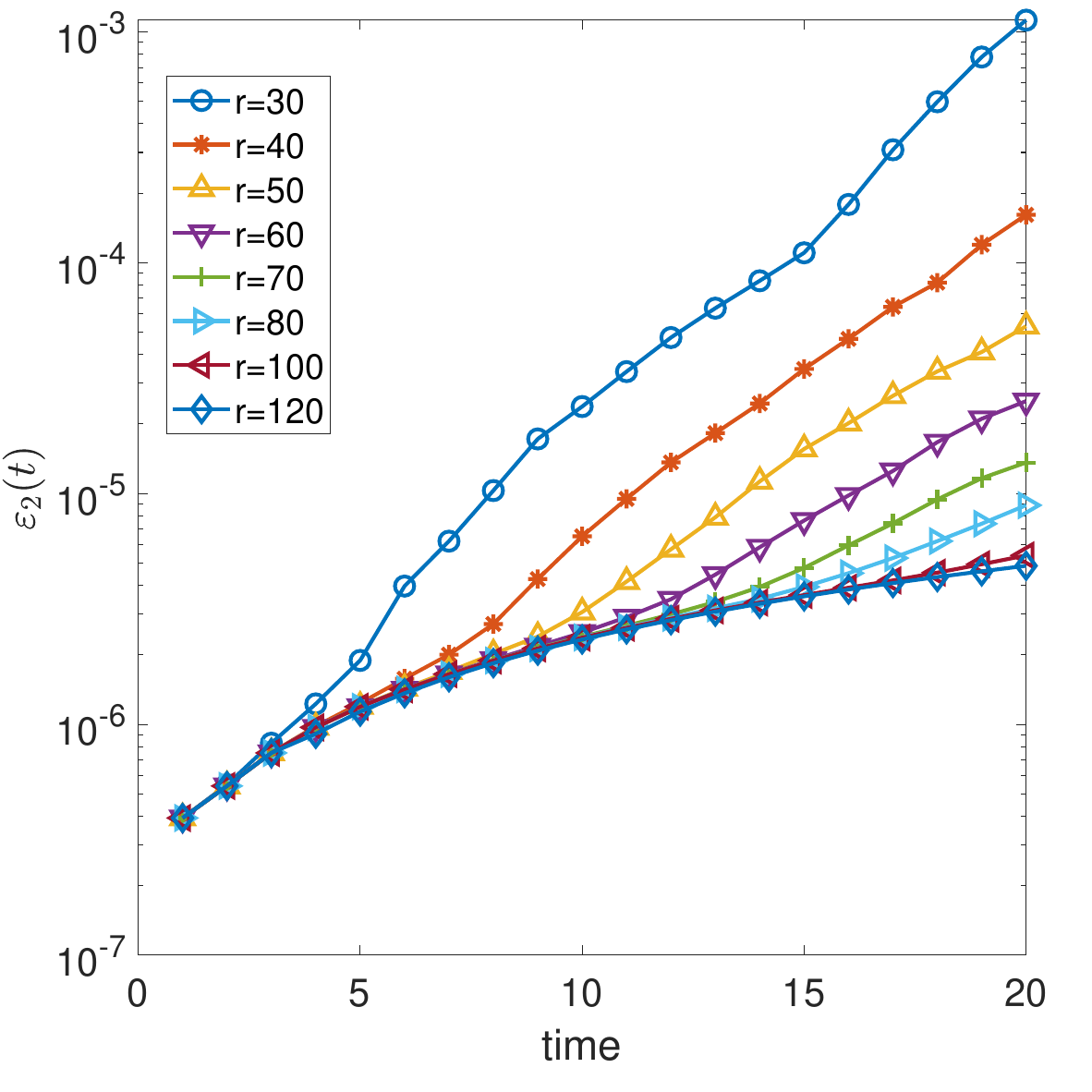}}
	\subfigure[$\varepsilon_{2}(t)$ under $r=100$ and different $N_t$. \label{fig::2Dbar_error_dt_e2}]
	{\includegraphics[width=0.45\textwidth,height=0.4\textwidth]{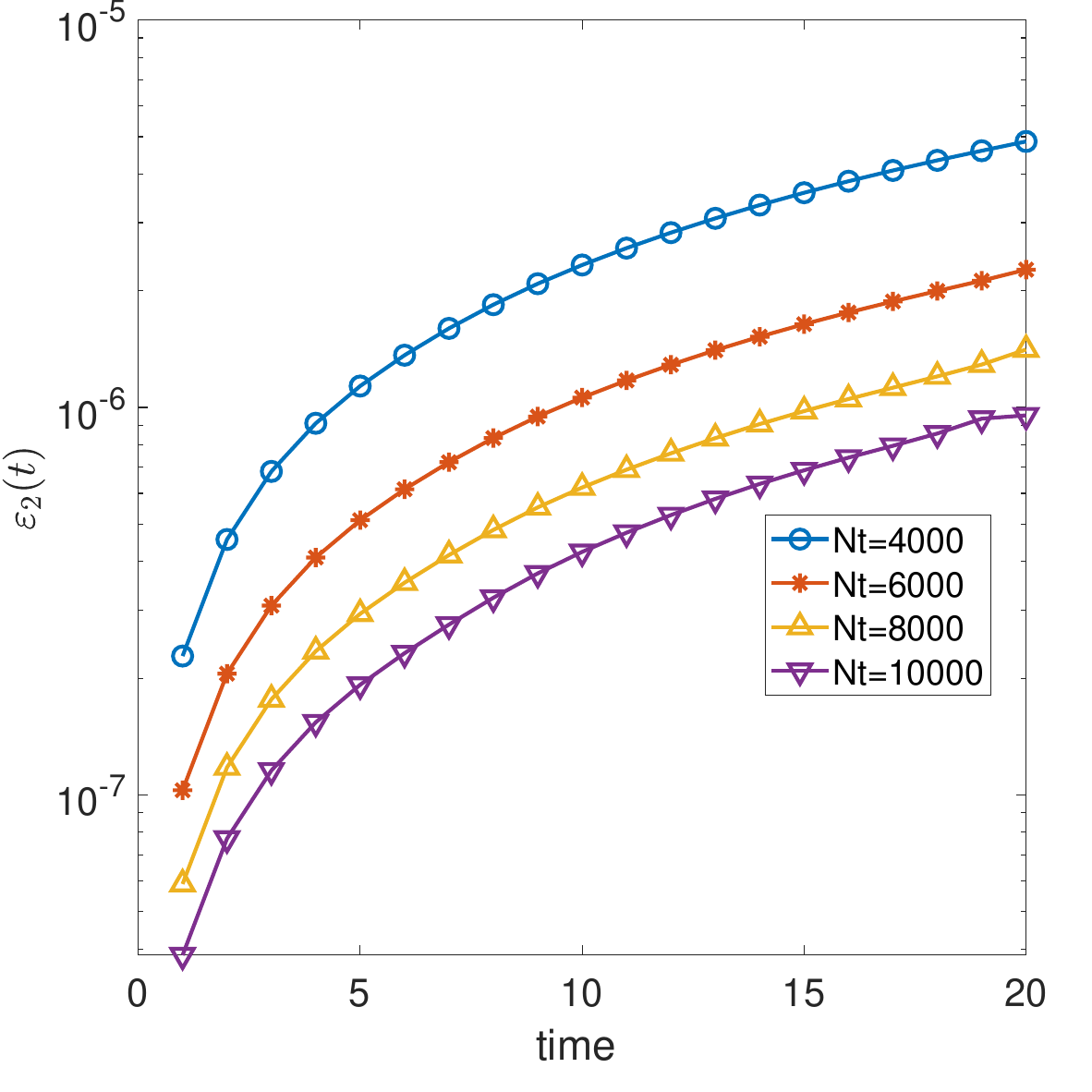}}
	\caption{2D Gaussian barrier scattering.}
	\label{fig::2Dbar_error_plot}
\end{figure}

A quantitative observation is conducted in \cref{fig::2Dbar_error_plot}. 
In this numerical example, the grid parameters are set as $N_x=200,N_k=400,$ and $R=80$ is estimated via SVD.
For reference solution, the same spatial grid subdivision parameters $N_x, N_k, N_y$ and the fixed time step $\Delta t=1/200$ are adopted as before. 
\Cref{fig::2Dbar_approxr} displays the $L_2$ error between the reference SEM solution and its best rank-$r$ approximation for $r=20,40,60,80,100,120,140.$
From this we can observe that the rank of this numerical test grows very rapidly over time, and relative to the grid size, the problem can be considered non‑low‑rank. 
For DLRA solution, we mainly focus on the convergence with respect to time step $\Delta t$ and rank $r$, several groups of under $\Delta t=1/200,1/300,1/400,1/500$ and $r=30,40,50,60,70,80,100,120$ are performed. 
\Cref{fig::2Dbar_error_r_e2} gives the relationship between error $\varepsilon_{2}(t)$ and rank $r$ under the same time step $\Delta t=1/500$. 
When the rank $r\le 100$, the error $\varepsilon_{2}(t)$ decreases noticeably as $r$ increase. 
However, the reduction of $\varepsilon_{2}(t)$ becomes much less pronounced or even ceases for $r>100$. 
This behavior occurs because the error is now dominated by the time step size $\Delta t$. 
We also examine the influence of the time step $\Delta t$ on the error, then \cref{fig::2Dbar_error_dt_e2} shows the dependence of the error $\varepsilon_{2}(t)$ on $\Delta t$ for a fixed rank $r=100$. 
We observe that errors accumulate progressively over time, and as the time step is reduced, $\varepsilon_{2}(t)$ decreases substantially while exhibiting a clear order of convergence.
The behavior of the other error indicator $\varepsilon_{\infty}(t)$ is similar to $\varepsilon_{2}(t)$. 
Moreover, the convergence rate with respect to $\Delta t$ presented in \cref{fig::2Dbar_rate}, clearly exhibits second‑order convergence as well.

\begin{figure}[h!] 
	\centering
	\subfigure[Initial value \cref{degenerate_initial_value}. \label{fig::4Drate_barde}]
	{\includegraphics[width=0.45\textwidth,height=0.4\textwidth]{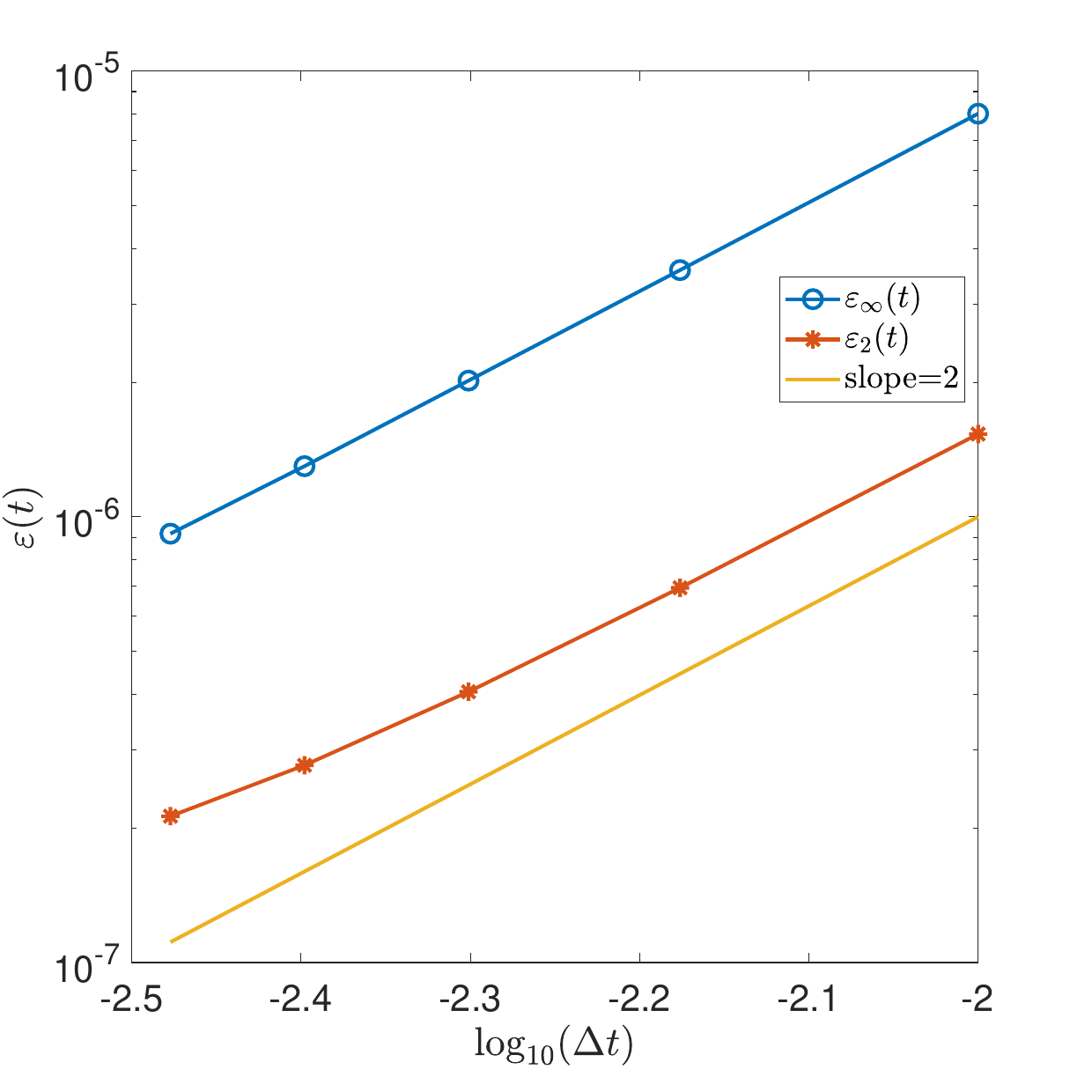}}
	\subfigure[Initial value \cref{common_initial_value}. \label{fig::4Drate_barfree}]
	{\includegraphics[width=0.45\textwidth,height=0.4\textwidth]{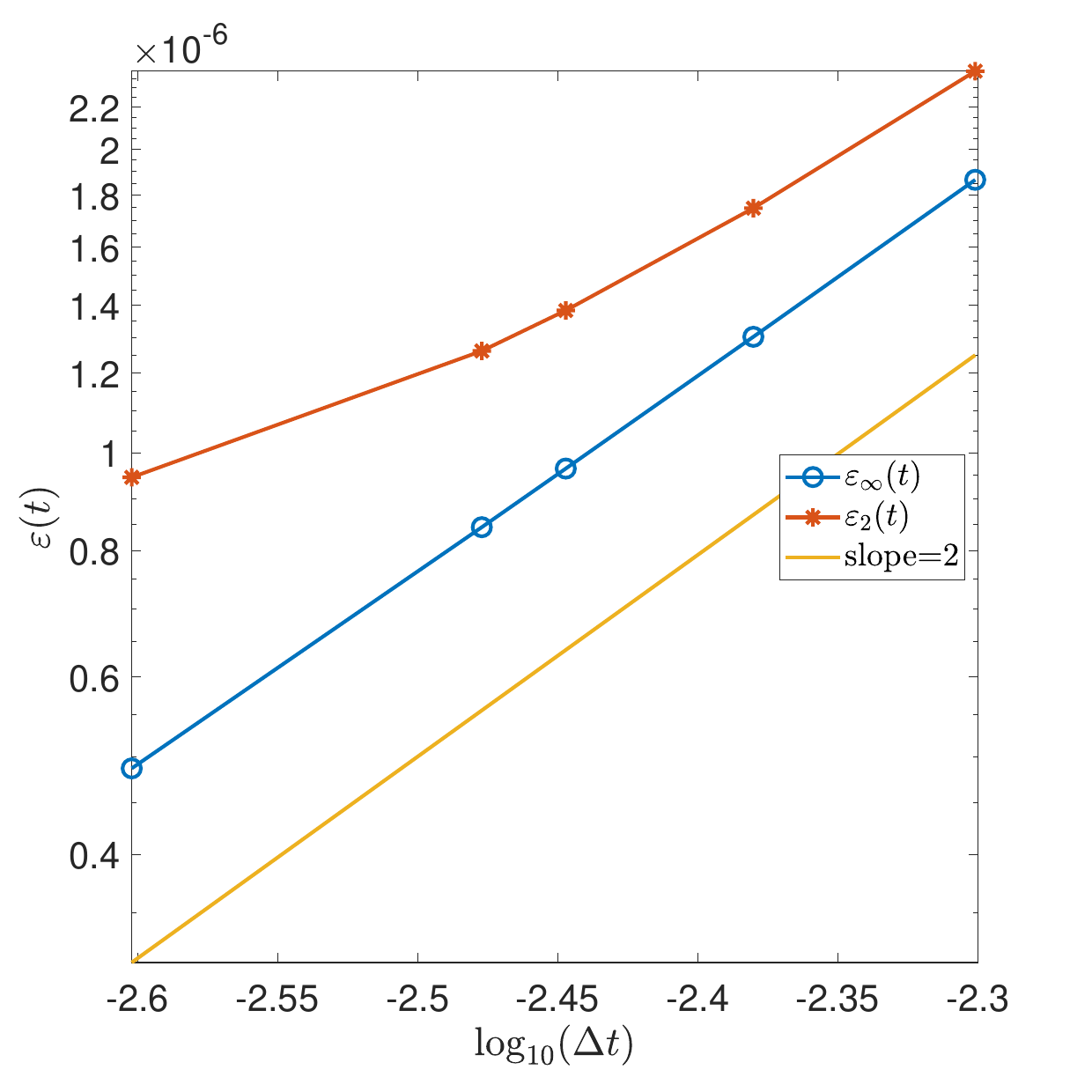}}
	\caption{4D Gaussian barrier scattering: The convergence order with respect to the time step $\Delta t$. The left figure shows the convergence order for the case with degeneracy in the $x_2$-$k_2$ direction, while the right figure displays the convergence order for the free advection in the $x_2$-$k_2$ direction. In both cases, second-order numerical accuracy can be observed.}
	\label{fig::4Drate_barfree_total}
\end{figure}

Such Gaussian barrier scattering can be extended to the 4D situation easily. 
The Wigner equation $f(x_1,x_2,k_1,k_2,t)$ with 2D potential $V(x_1,x_2)=V_1(x_1)+V_2(x_2)$ has the form 
\begin{align}
	f(x_1,x_2,k_1,k_2,t)=f_1(x_1,k_1,t)f_2(x_2,k_2,t),
\end{align}
where $f_i(x_i,k_i,t)$ is 2D Wigner function corresponding to the 1D potential $V_i(x_i)$ for $i=1,2.$ 
For instance, we adopt computational domain $\Omega=\Omega_{\bx}\times \Omega_{\bk}=[0,60]^2\times [-4,4]^2$ with spacing $\Delta y_1=\Delta y_2= \pi/4$.
The potential function is chosen as
\begin{align}
	V(x_1,x_2)=H_P\exp\left[-\frac{(x_1-x_{1,B})^2}{2w^2}\right], \label{4Dbarfree_potent}
\end{align}
where $H_P=0.3, x_{1,B}=30$ and $w=1$. The initial Wigner function is set as degenerate function
\begin{align}
	f(x_1,x_2,k_1,k_2,0)=\frac{1}{\pi}\frac{1}{L_{x_2}}\frac{1}{L_{k_2}}\exp\left[-\frac{(x_1-x_{1,0})^2}{2a_1^2}\right]\exp\left[-2a_1^2(k_1-k_{1,0})^2\right], \label{degenerate_initial_value}
\end{align}
where $x_{1,0}=15, k_{1,0}=0.7, a_1=2.825$, $L_{x_2}=|\Omega_{x_2}|=60$ and $L_{k_2}=|\Omega_{k_2}|=4$. 
In this case, the Wigner function $f(x_1,x_2,k_1,k_2,t)$ is degenerate in $x_2$-direction and $k_2$-direction, and the reduced Wigner function $W_1(x,k,t)$ is the same as that of the previous numerical test. 
The remaining parameters are the same as the previous 2D case. 
To study the convergence rate with respect to time step $\Delta t$, we fix $N_{x_1}=N_{x_2}=N_{k_1}=N_{k_2}=100$ and $R=63$ is estimated by random SVD. 
The order of convergence is presented in \cref{fig::4Drate_barde},  we can clearly observe the second-order convergence order as well. 

Similarly, if we choose the same potential function in \cref{4Dbarfree_potent} with initial Wigner function as
\begin{equation}
	\begin{aligned}
		f(x_1,x_2,k_1,k_2,0)=\frac{1}{\pi^2}&\exp\left[-\frac{(x_1-x_{1,0})^2}{2a_1^2}-\frac{(x_2-x_{2,0})^2}{2a_2^2}\right]\\
		\times&\exp\left[-2a_1^2(k_1-k_{1,0})^2-2a_2^2(k_2-k_{2,0})^2\right], 
	\end{aligned}\label{common_initial_value}
\end{equation}
where $x_{1,0}, k_{1,0}$ and $a_1$ are in the same way as before, $x_{2,0}=45,k_{2,0}=-0.7$ and $a_2=2.825$.
In this case, the reduced Wigner function $W_1(x,k,t)$ is the same as that of the previous numerical test and $W_2(x,k,t)$ evolves as the free convection of an initial GWP over time.
The remaining parameters are the same as the previous case except the rank $r$ is chosen as $r=400$, which is much larger than the previous.
The convergence order, shown in \cref{fig::4Drate_barfree}, also clearly exhibits second‑order accuracy.

\subsection{Electron-electron scattering}
Now we turn to discuss a more challenging problem. 
Consider that two electrons interacting through the repulsive Coulomb force, the potential function is given by the soft Coulomb potential
\begin{align}
	V_{s}(x_1,x_2)=\frac{1}{\sqrt{(x_1-x_2)^2+\epsilon}}, \label{soft Coulomb}
\end{align}
where $\epsilon=1$ is the soft parameter.

\begin{figure}[h!]
	\centering
	\subfigure[Initial value: two Fermions\label{fig::Fermion}]
	{\includegraphics[width=0.48\textwidth,height=0.4\textwidth]{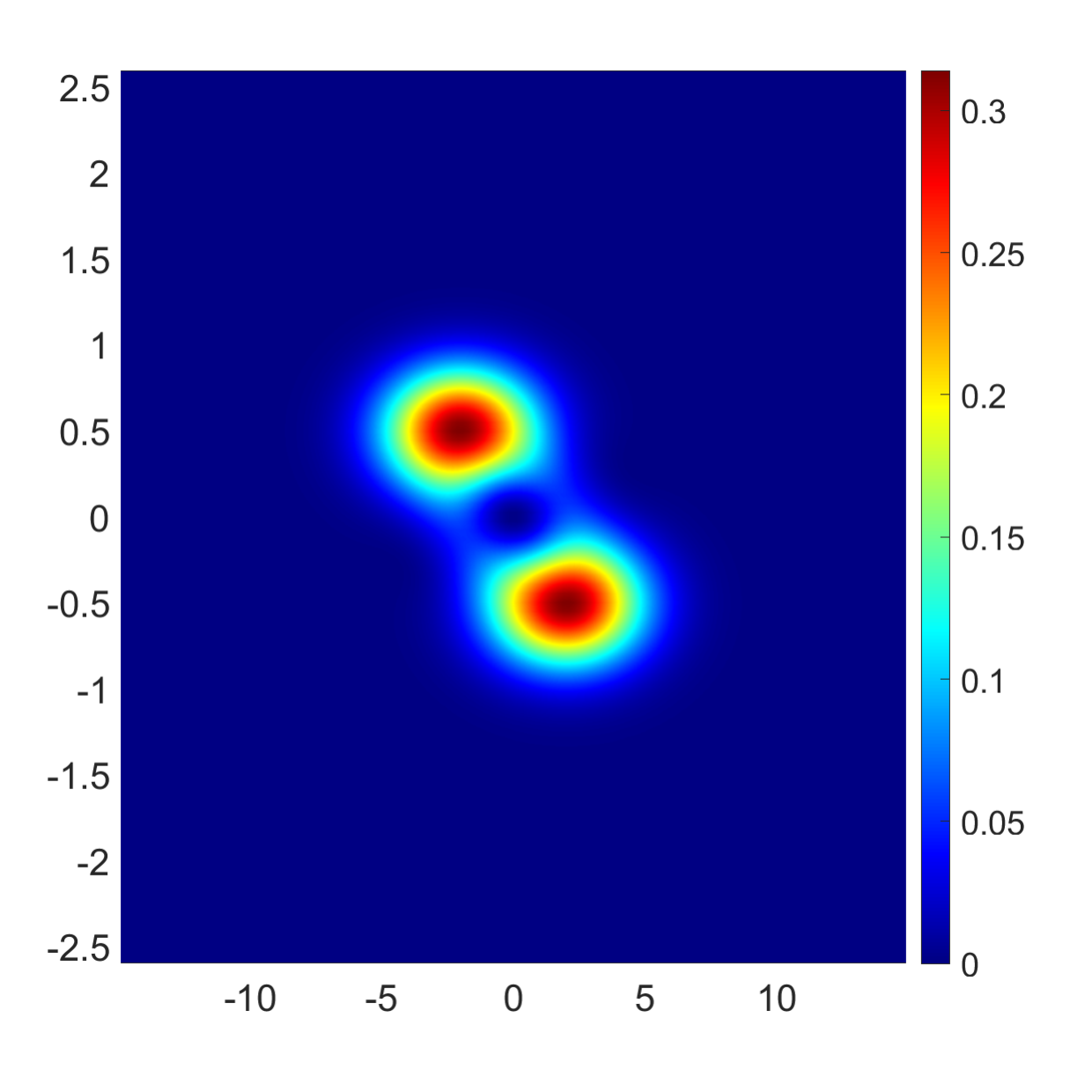}}
	\subfigure[Electron-electron scattering: The convergence order with respect to the time step $\Delta t$. \label{fig::4Dee_rate}]
	{\includegraphics[width=0.48\textwidth,height=0.4\textwidth]{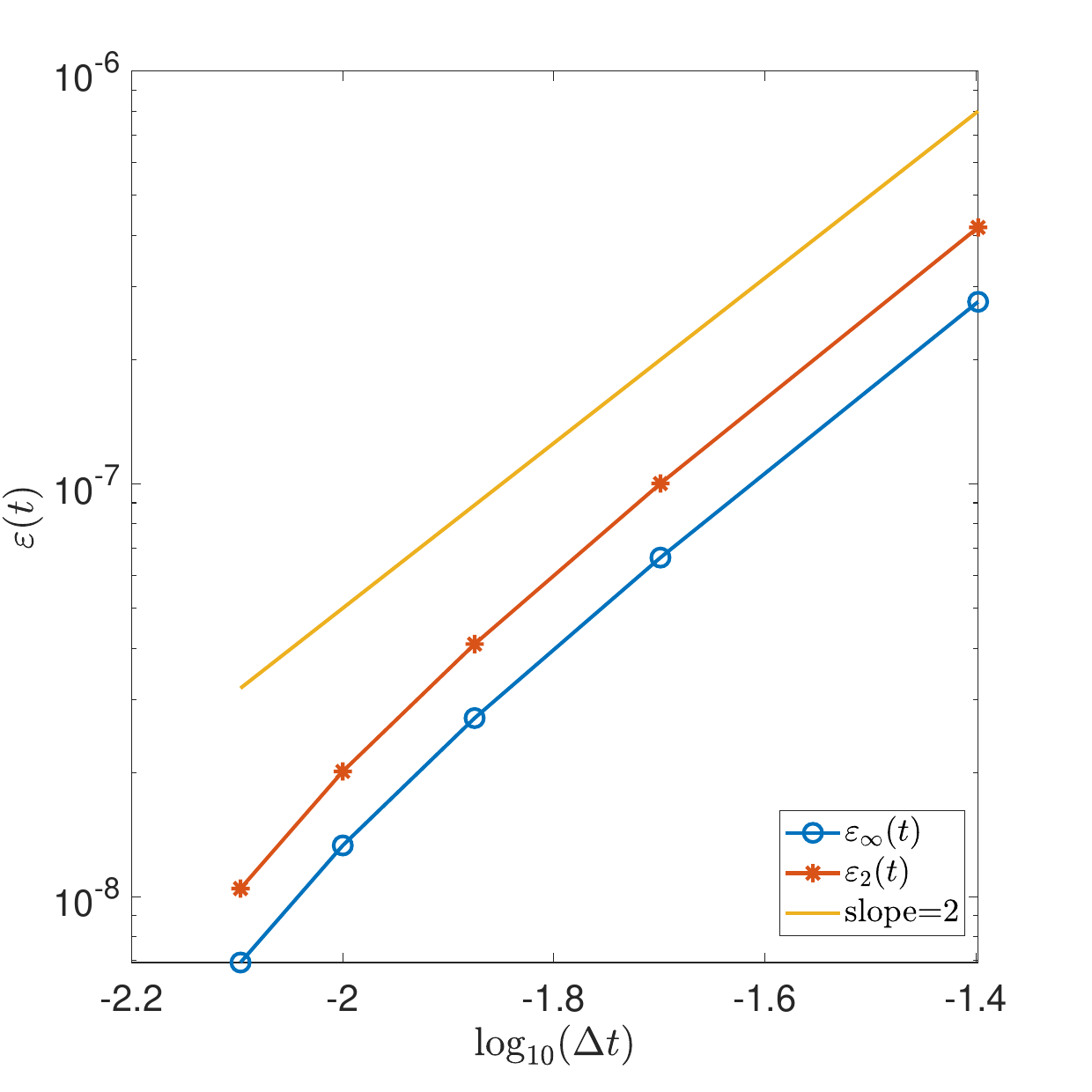}}
	\caption{The reduced Wigner function $W_1(x,k,t)+W_2(x,k,t)$ for two fermions and the convergence order.}
	\label{fig::Fermion_and_rate}
\end{figure}

The initial value of the Wigner function to simulate two electronic states is chosen as a Fermion system
\begin{equation}
	\begin{aligned}
	& f_0^{\text {Fermion }}\left(x_1, x_2, k_1, k_2\right) \\
	= & \frac{1}{2 \pi^2} \exp \left[-\frac{\left(x_1-x_1^0\right)^2}{2 a^2}-\frac{\left(x_2-x_2^0\right)^2}{2 a^2}-2 a^2\left(k_1-k_1^0\right)^2-2 a^2\left(k_2-k_2^0\right)^2\right] \\
	& +  \frac{1}{2 \pi^2} \exp \left[-\frac{\left(x_1-x_2^0\right)^2}{2 a^2}-\frac{\left(x_2-x_1^0\right)^2}{2 a^2}-2 a^2\left(k_1-k_2^0\right)^2-2 a^2\left(k_2-k_1^0\right)^2\right] \\
	& - \frac{1}{\pi^2} \exp \left[-\frac{\left(x_1-x_1^0\right)^2+\left(x_1-x_2^0\right)^2+\left(x_2-x_1^0\right)^2+\left(x_2-x_2^0\right)^2}{4 a^2}\right] \\
	& \times  \exp \left[\frac{\left(x_1^0-x_2^0\right)^2}{4 a^2}-2 a^2\left(k_1-\frac{k_1^0+k_2^0}{2}\right)^2-2 a^2\left(k_2-\frac{k_1^0+k_2^0}{2}\right)^2\right] \\
	& \times  \cos \left[\left(x_1^0-x_2^0\right)\left(k_1-k_2\right)-\left(k_1^0-k_2^0\right)\left(x_1-x_2\right)\right].
	\end{aligned} 
	\label{Fermion}
\end{equation}
In the following computations, the Fermion system is taken as the initial condition, with the parameters set as $x_1^0=-2,x_2^0=2,k_1^0=0.5,k_2^0=-0.5, a=2,$ as illustrated in \cref{fig::Fermion}.

The first run is to demonstrate the numerical convergence of the proposed DLRA method by setting computational domain $\Omega=\Omega_{\bx}\times\Omega_{\bk}=[-15,15]^2 \times [-5\pi/6,5\pi/6]^2$ with spacings $\Delta y_1=\Delta y_2=1.2$. 
The number of grid points in each spatial direction is set to $N_{x_1}=N_{x_2}=150,$ and in each momentum direction to $ N_{k_1}=N_{k_2}=100$ with final time $T=4$. 
The rank $r=200$ is chosen sufficiently large to reveal the convergence order with respect to the time step $\Delta t$, while $R=58$ is estimated by random SVD. 
The reference solution is computed using advective-spectral-mixed (ASM) scheme \cite{xiong2016advectivespectralmixed} with the same parameters listed above. 
\Cref{fig::4Dee_rate} presents the results, clearly showing that our DLRA algorithm achieves second‑order convergence.
However, choosing an excessively large rank $r$ leads to high computational cost, especially in the $S$-step, which contradicts the original purpose of adopting a low‑rank approach. 
In this numerical test, a rank of $r=150$ is sufficient to produce an approximate solution that is visually indistinguishable from the reference.

\begin{figure}[h]
	\centering
	\subfigure[t=1.\label{fig::t=1}]
	{\includegraphics[width=0.48\textwidth,height=0.35\textwidth]{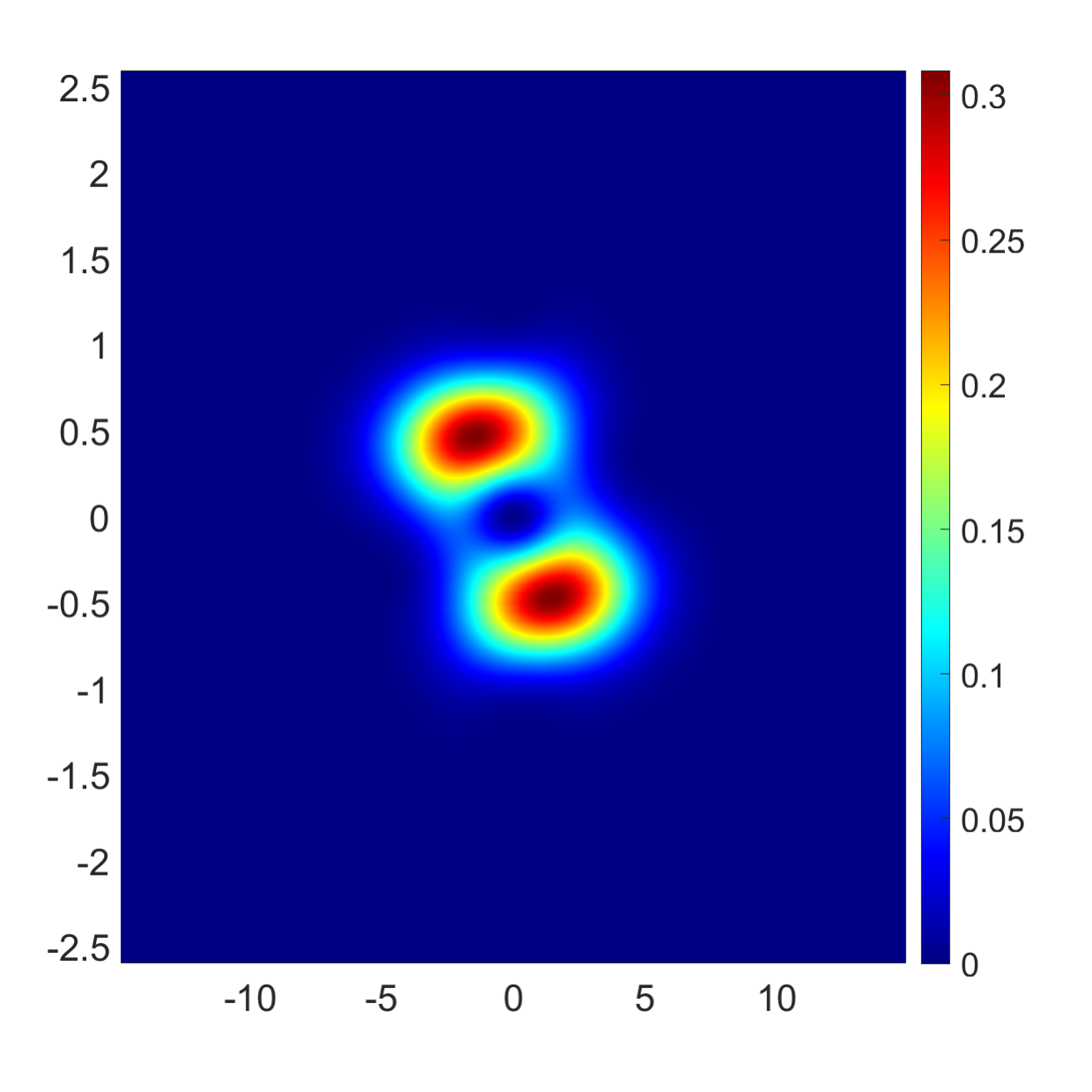}}
	\subfigure[t=2.\label{fig::t=2}]
	{\includegraphics[width=0.48\textwidth,height=0.35\textwidth]{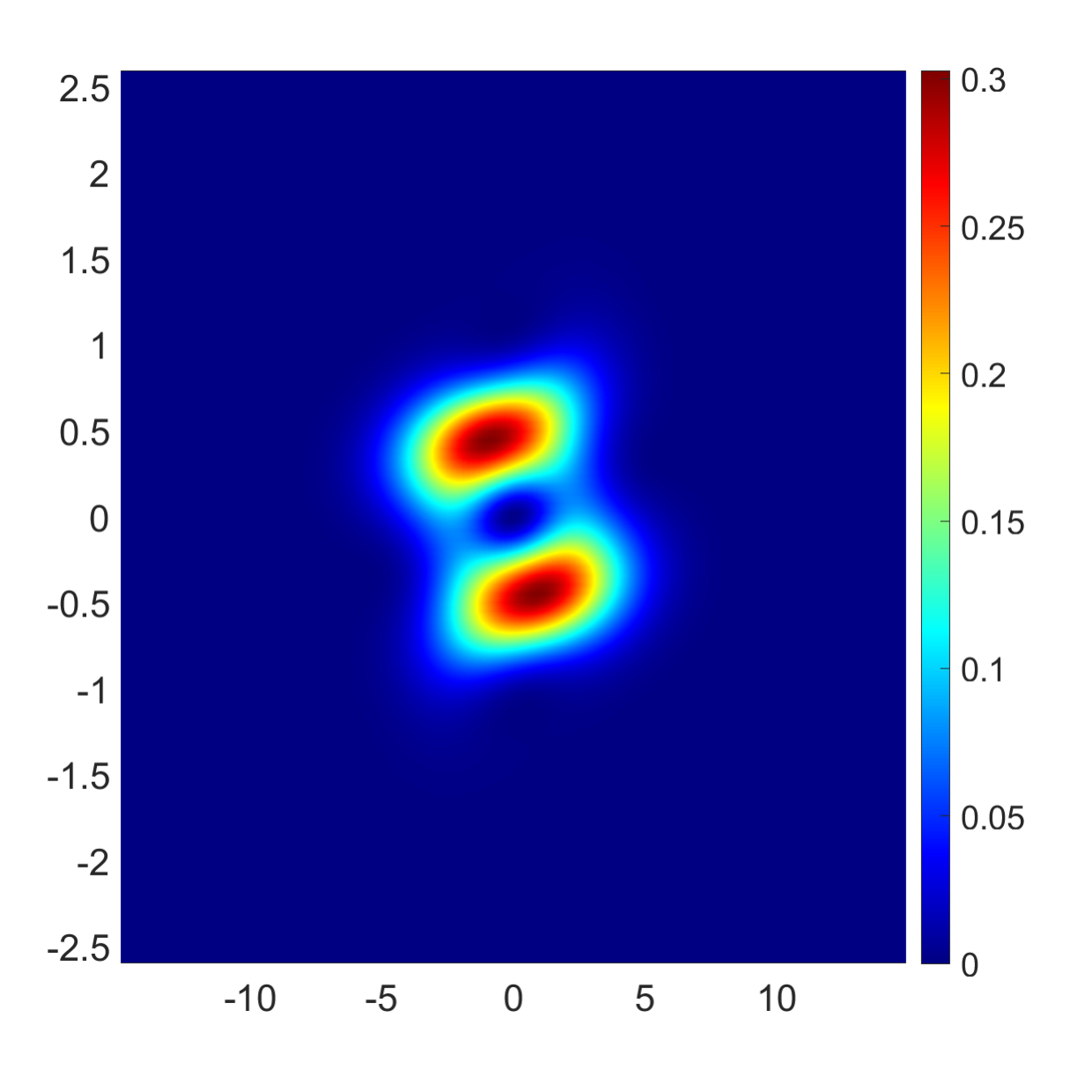}}\\
	\subfigure[t=3.\label{fig::t=3}]
	{\includegraphics[width=0.48\textwidth,height=0.35\textwidth]{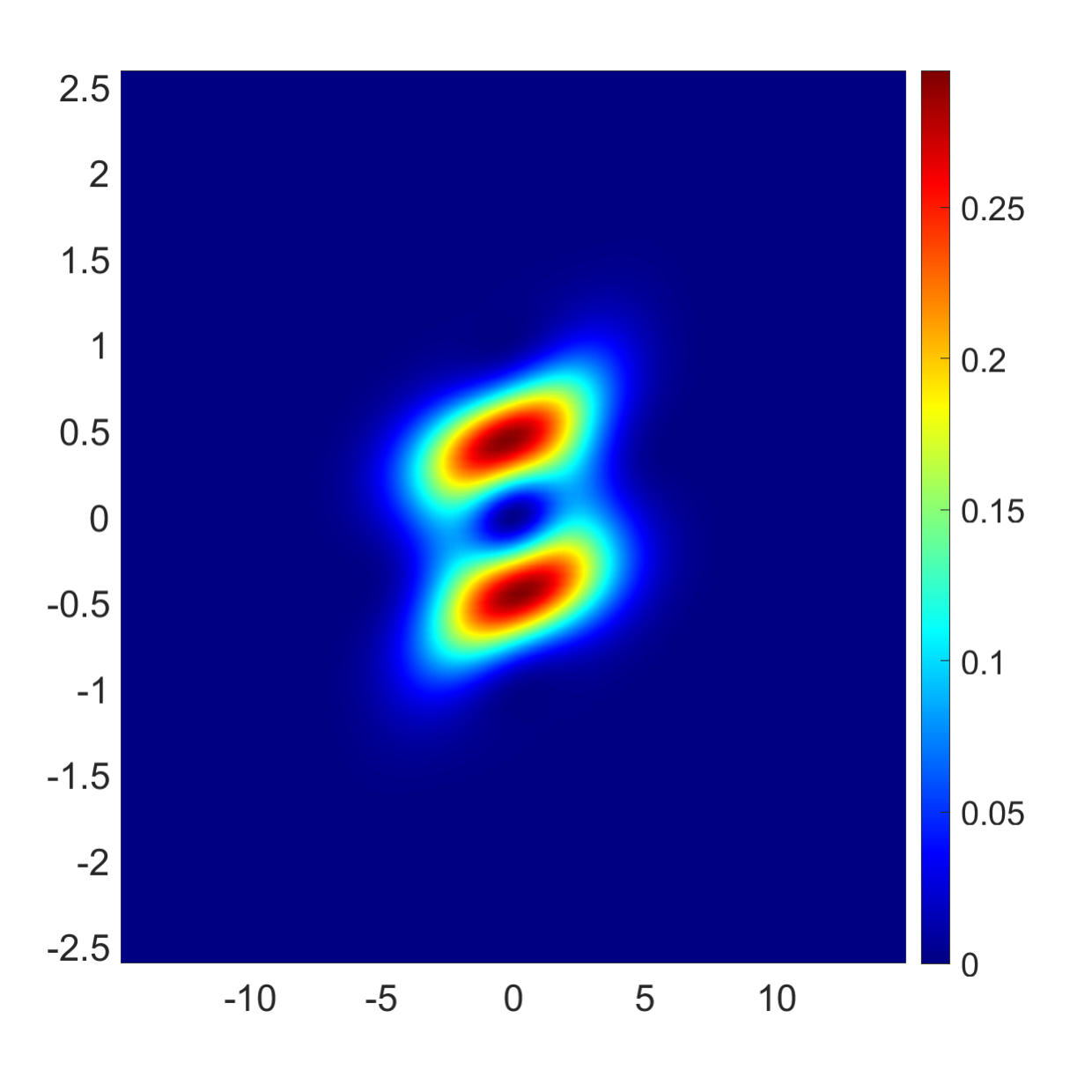}}
	\subfigure[t=4.\label{fig::t=4}]
	{\includegraphics[width=0.48\textwidth,height=0.35\textwidth]{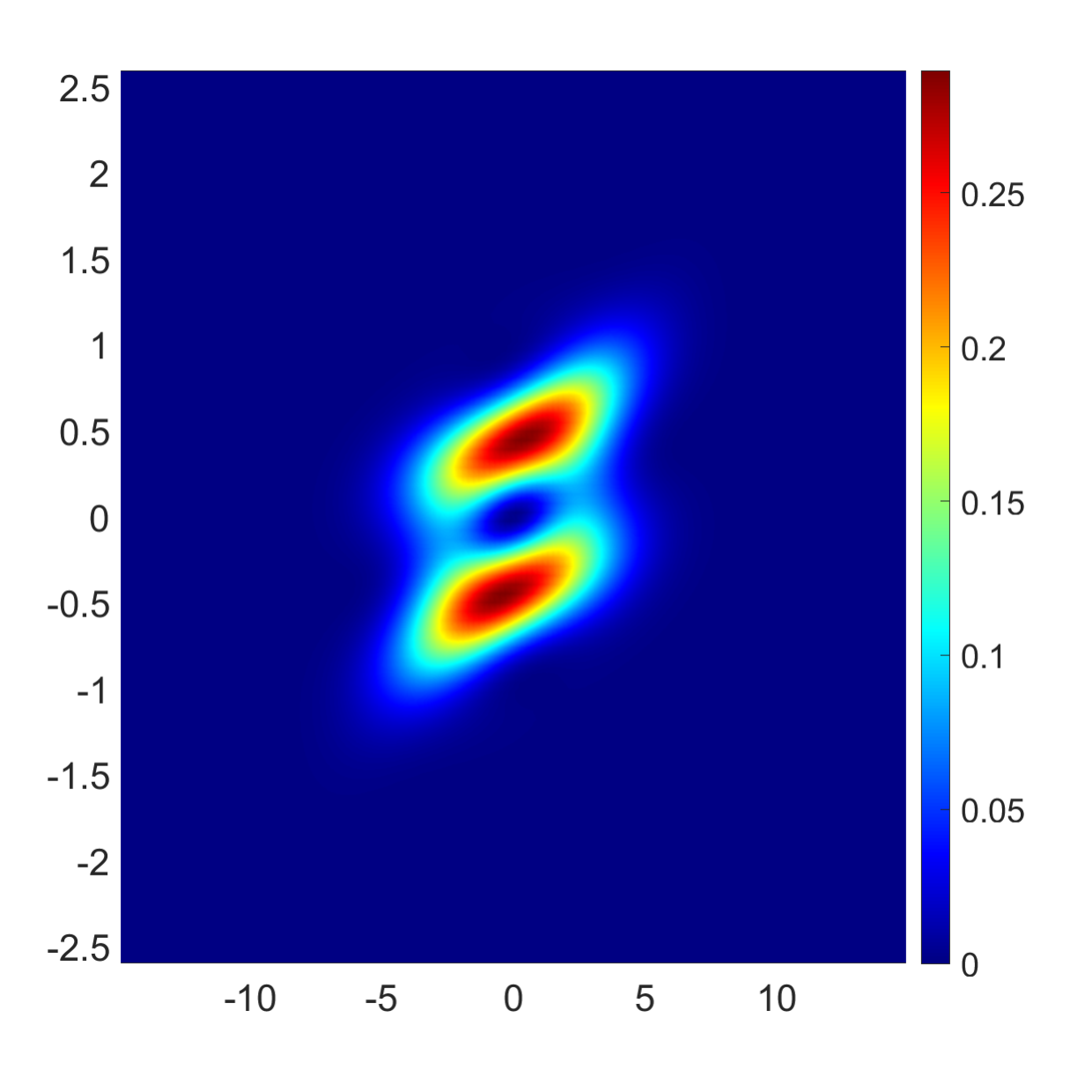}}\\
	\subfigure[t=5.\label{fig::t=5}]
	{\includegraphics[width=0.48\textwidth,height=0.35\textwidth]{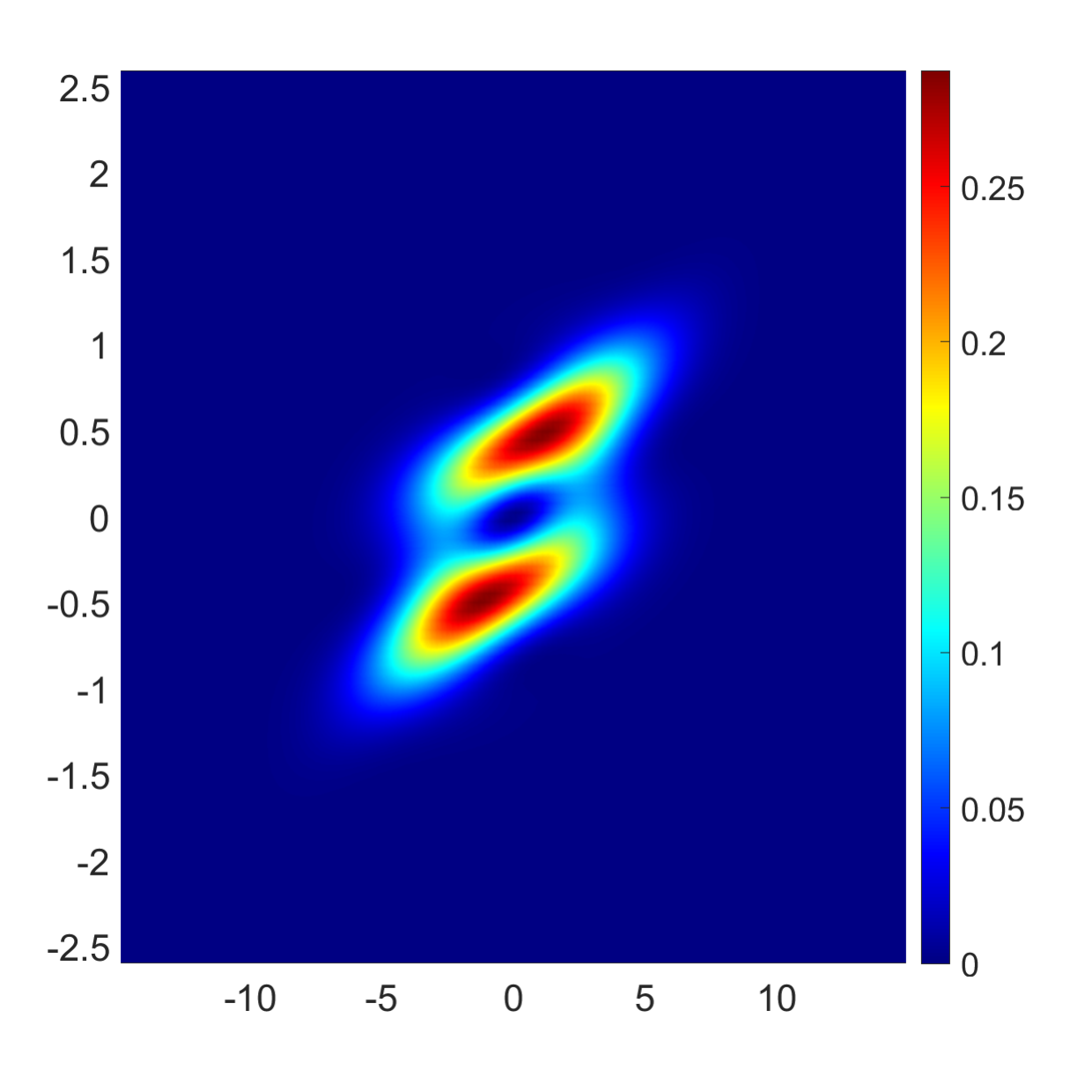}}
	\subfigure[t=6.\label{fig::t=6}]
	{\includegraphics[width=0.48\textwidth,height=0.35\textwidth]{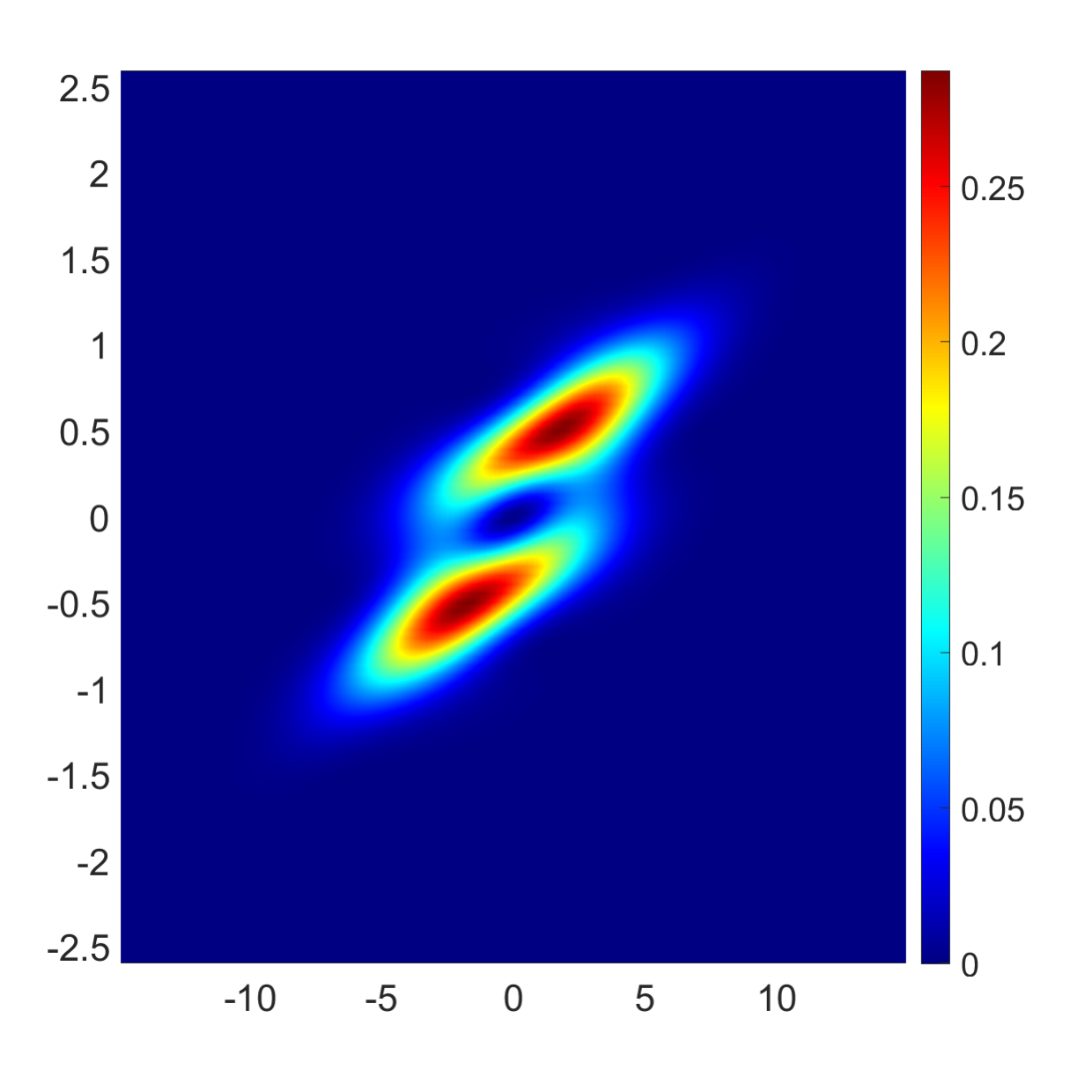}}
	\caption{Electron-electron scattering: $W_1(x,k,t)+W_2(x,k,t)$ at different time.}
	\label{fig::ee}
\end{figure}

Then, \cref{fig::ee} shows the evolution of the reduced Wigner function $W_1(x,k,t)+W_2(x,k,t)$ under the same parameters listed above, up to final time $T=6$. 
The reduced Wigner function of the two electrons moving towards each other and decelerate before $t=4$ is suppressed in the $|k|\le 1$ region and has greater expansion in the $x$-space due to Coulomb potential. Since then, two electrons tend to scattered and separated from each other due to repulsive interaction, which is line with our expectations. 
Moreover, it should be emphasized that at $t=0$, a hole structure known as the Fermi hole, appears at the center $\left((x_1^0+x_2^0)/2, (k_1^0+k_2^0)/2\right)$ between the two electrons as a direct manifestation of the Pauli exclusion principle.
Furthermore, this Fermi hole persists throughout the evolution due to the strong correlation between the electrons. 

Finally, we employ this example to illustrate the advantages of the DLRA method in handling such high‑dimensional problems.
\Cref{reality compare} compares the theoretical memory requirement and the running time per step for this problem under different ranks $r$, using the same parameters as above, for both the reference ASM scheme and our DLRA method.
It can be observed that the DLRA method outperforms the reference full‑grid method in both theoretical memory requirement and computational time.
For example, with $r=150$ being sufficient for this problem, the memory required by the DLRA algorithm is two orders of magnitude smaller, and the running time is one order of magnitude faster.

\begin{table}[h]
	\centering
	\renewcommand\arraystretch{1.2}
	\caption{Comparison of computational cost for electron–electron scattering simulations: full-grid (ASM) vs. DLRA with different ranks $r$.}
	\label{reality compare}
	\begin{tabular}{ >{\centering\arraybackslash}p{9em} | *{5}{>{\centering\arraybackslash}c} }
		\hline
		& Full grids (ASM) & $r=50$  & $r=100$ & $r=150$ & $r=200$ \\ [1pt]\hline
		\makecell[cc]{Theoretical memory \\ requirement (GB)} 
		& 23.45 & 0.04 & 0.07 & 0.10 & 0.13   \\ [1pt]\hline
		\makecell[cc]{Running time \\ per step (s)} 
		& 75.76 & 2.58 & 5.34 & 8.35 & 11.34   \\ \hline
	\end{tabular}
\end{table}

\subsection{A Helium-like system}
As the final example, we consider a Helium-like system composed of two electrons. 
Besides the repulsive Coulomb force, they are both attracted by a Helium atom at $\bx=\textbf{0}$. 
To describe the nucleon-electron interaction, we adopt the attractive soft-Coulomb potential \cite{lein2002strongfield}
\begin{align}
	V_{he}(x_1,x_2)=-\frac{Z}{\sqrt{|x_1|^2+\epsilon}}-\frac{Z}{\sqrt{|x_2|^2+\epsilon}},
\end{align}
with the atomic number $Z=2$ and the soft parameter $\epsilon=1.$ 
The initial value is taken to be the same Fermion system as in the previous case, using the identical parameter set, shown in \cref{fig::Fermion}.
And the other parameters are: the computational domain $\Omega=\Omega_{\bx}\times\Omega_{\bk}=[-15,15]^2\times[-5\pi/3,5\pi/3]^2$ with spacings $\Delta y_1=\Delta y_2=0.6, $
the number of grid points in each spatial direction is set to $N_{x_1}=N_{x_2}=150,$ and in each momentum direction to $ N_{k_1}=N_{k_2}=100,$ with $R=106$ is estimated by random SVD, while $r=400$ is chosen sufficiently large. 

\begin{figure}[h!]
	\centering
	\includegraphics[width=0.47\linewidth]{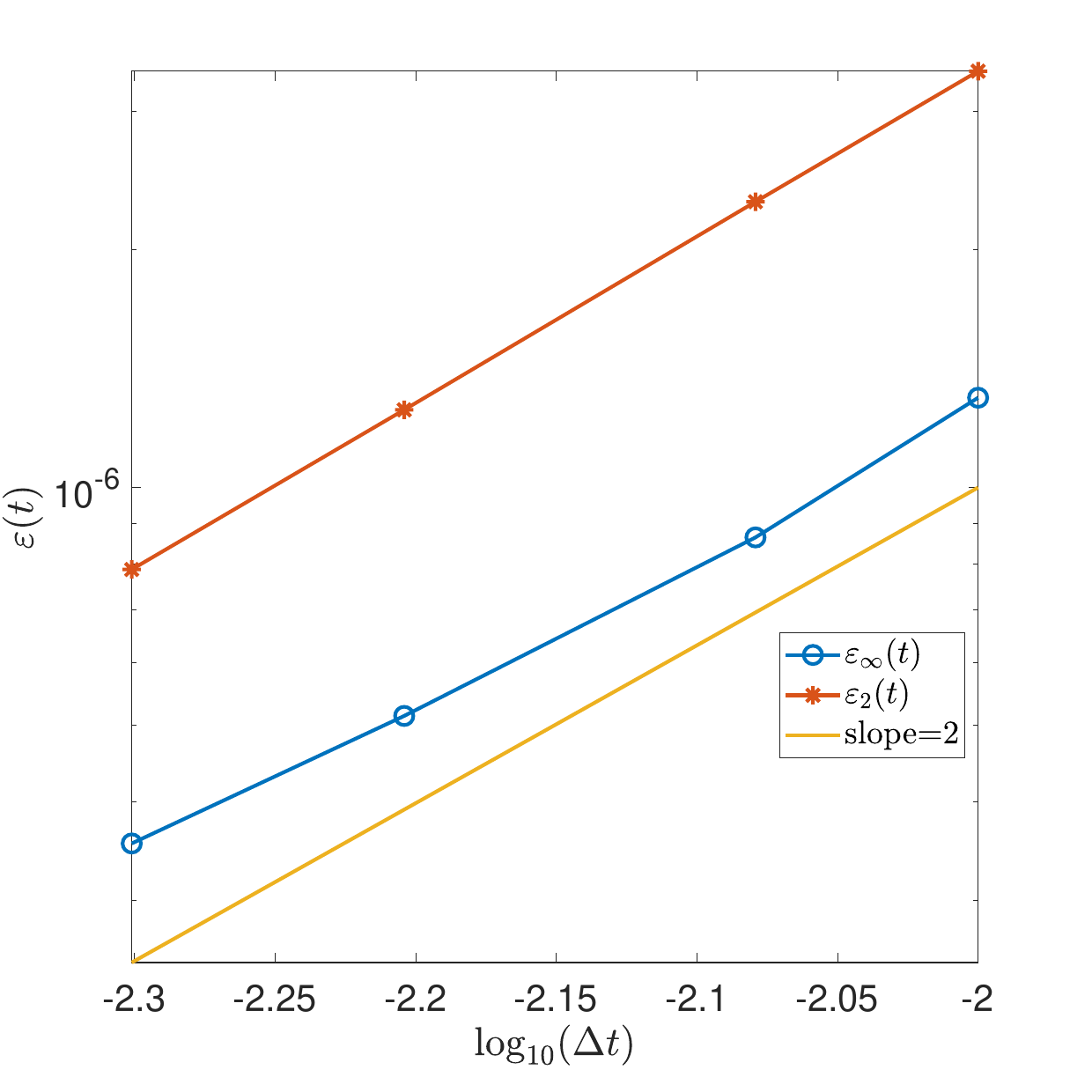}
	\caption{The reduced Wigner function $W_1(x,k,t)+W_2(x,k,t)$ for two fermions and the convergence order.}
	\label{fig::he_rate}
\end{figure}

The first run is carried out to demonstrate the numerical convergence of the proposed DLRA method. 
\Cref{fig::he_rate} presents the results, clearly showing that our DLRA method achieves second‑order convergence at time $T=1$. 
Similarly, in this test case, a sufficiently large rank is required to reliably measure the convergence order; rank $r=150$ is sufficient to compute an approximate solution that captures the relevant physical properties.

\begin{figure}[h!]
	\centering
	\subfigure[t=1.]
	{\includegraphics[width=0.32\textwidth,height=0.25\textwidth]{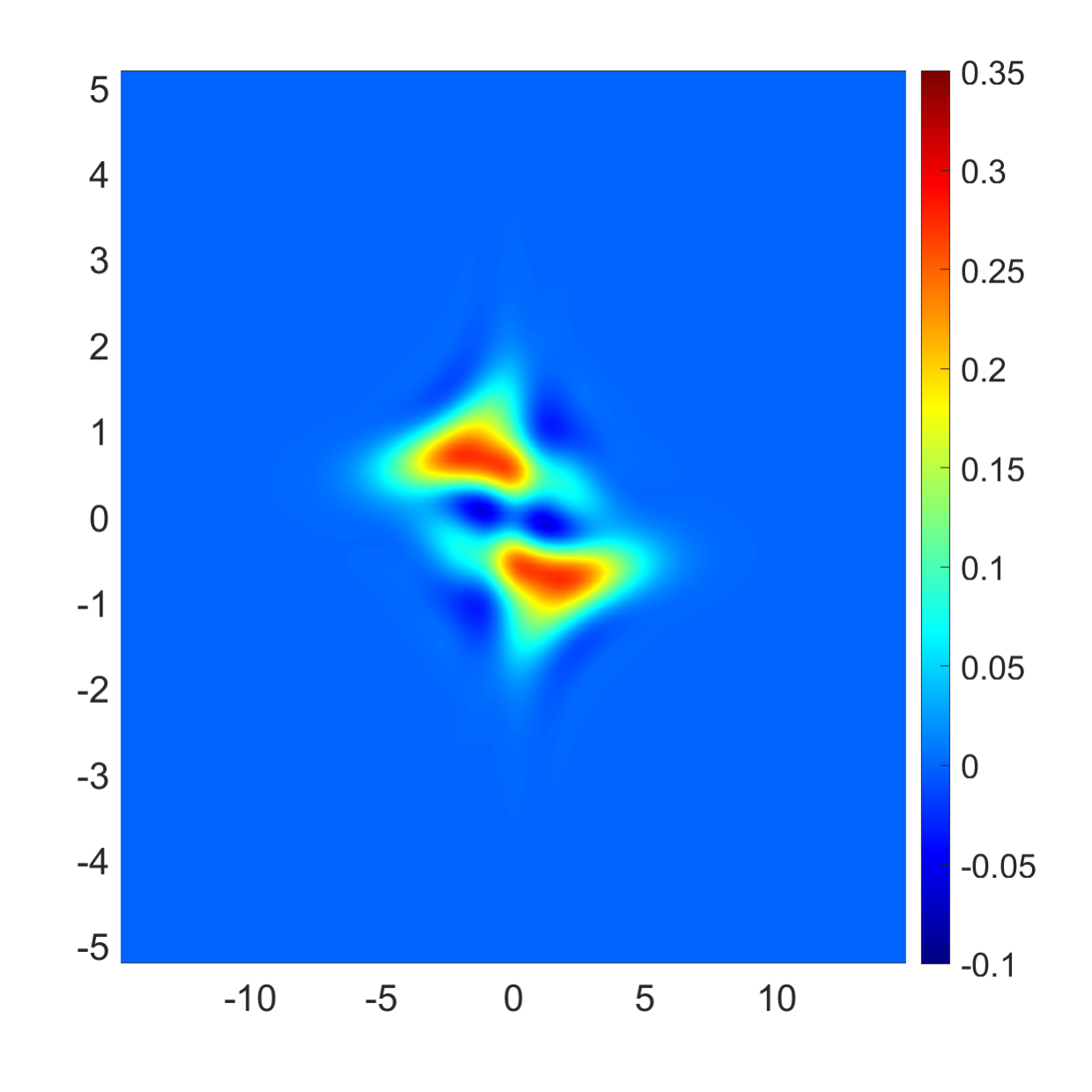}}
	\subfigure[t=2.]
	{\includegraphics[width=0.32\textwidth,height=0.25\textwidth]{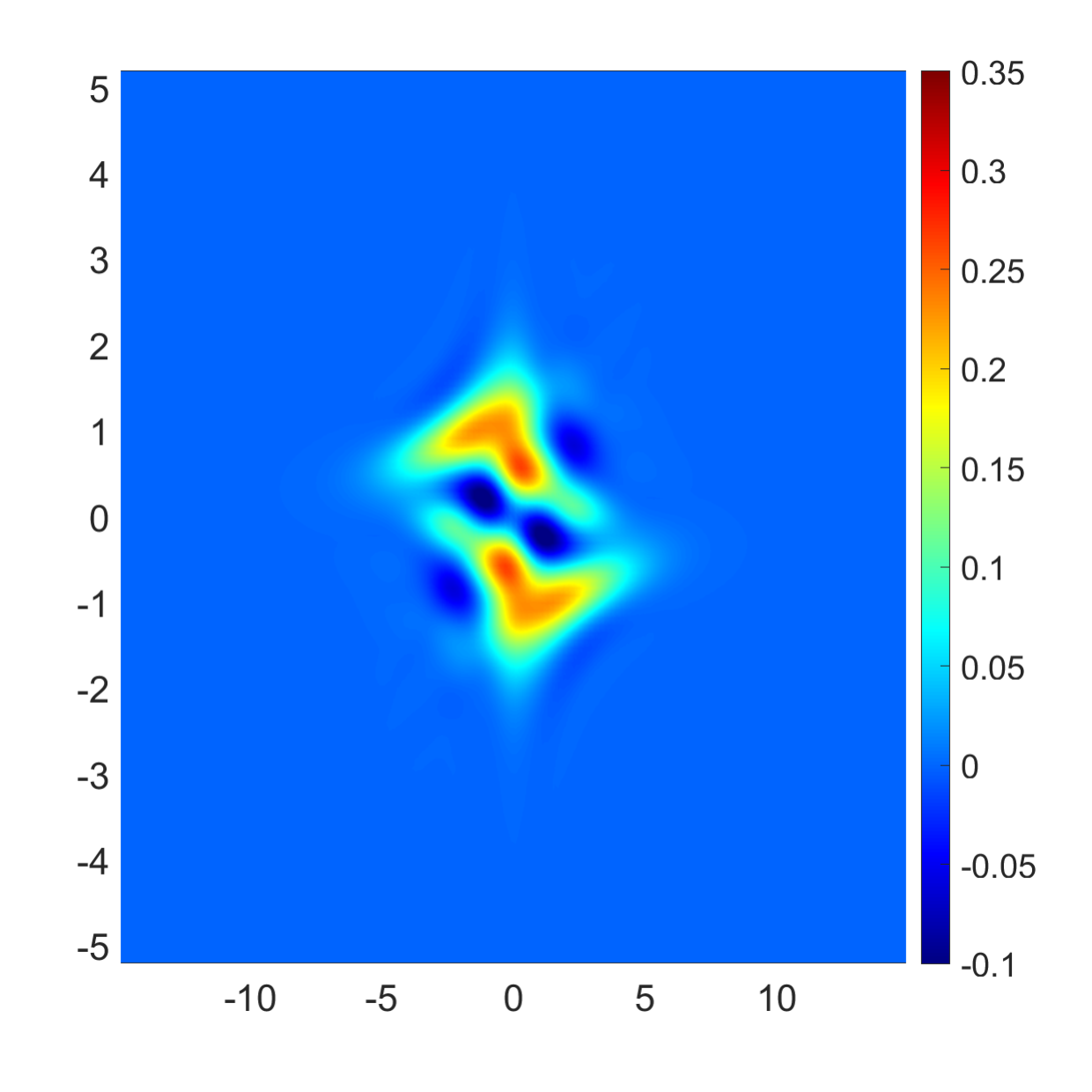}}
	\subfigure[t=3.]
	{\includegraphics[width=0.32\textwidth,height=0.25\textwidth]{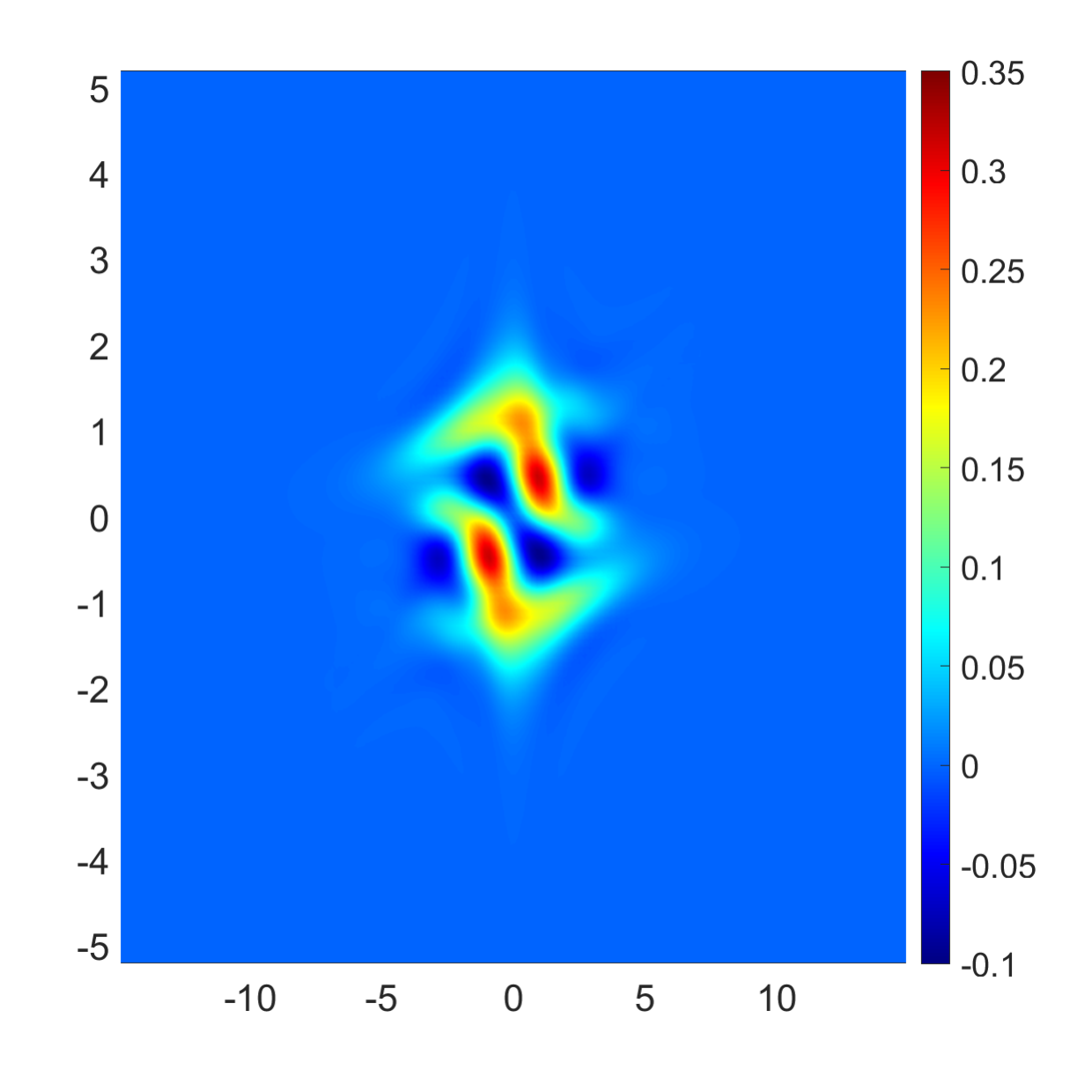}}\\
	\subfigure[t=4.]
	{\includegraphics[width=0.32\textwidth,height=0.25\textwidth]{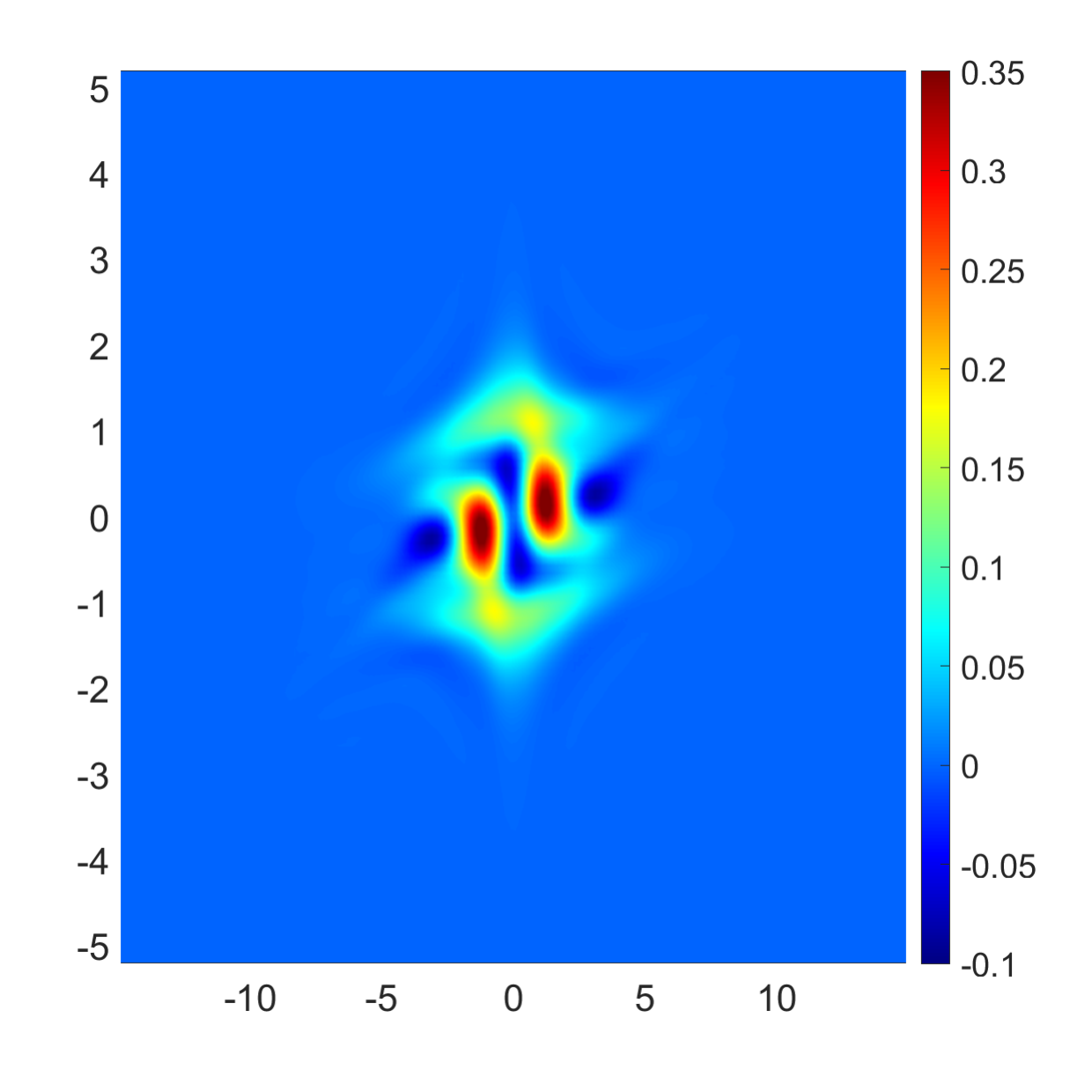}}
	\subfigure[t=5.]
	{\includegraphics[width=0.32\textwidth,height=0.25\textwidth]{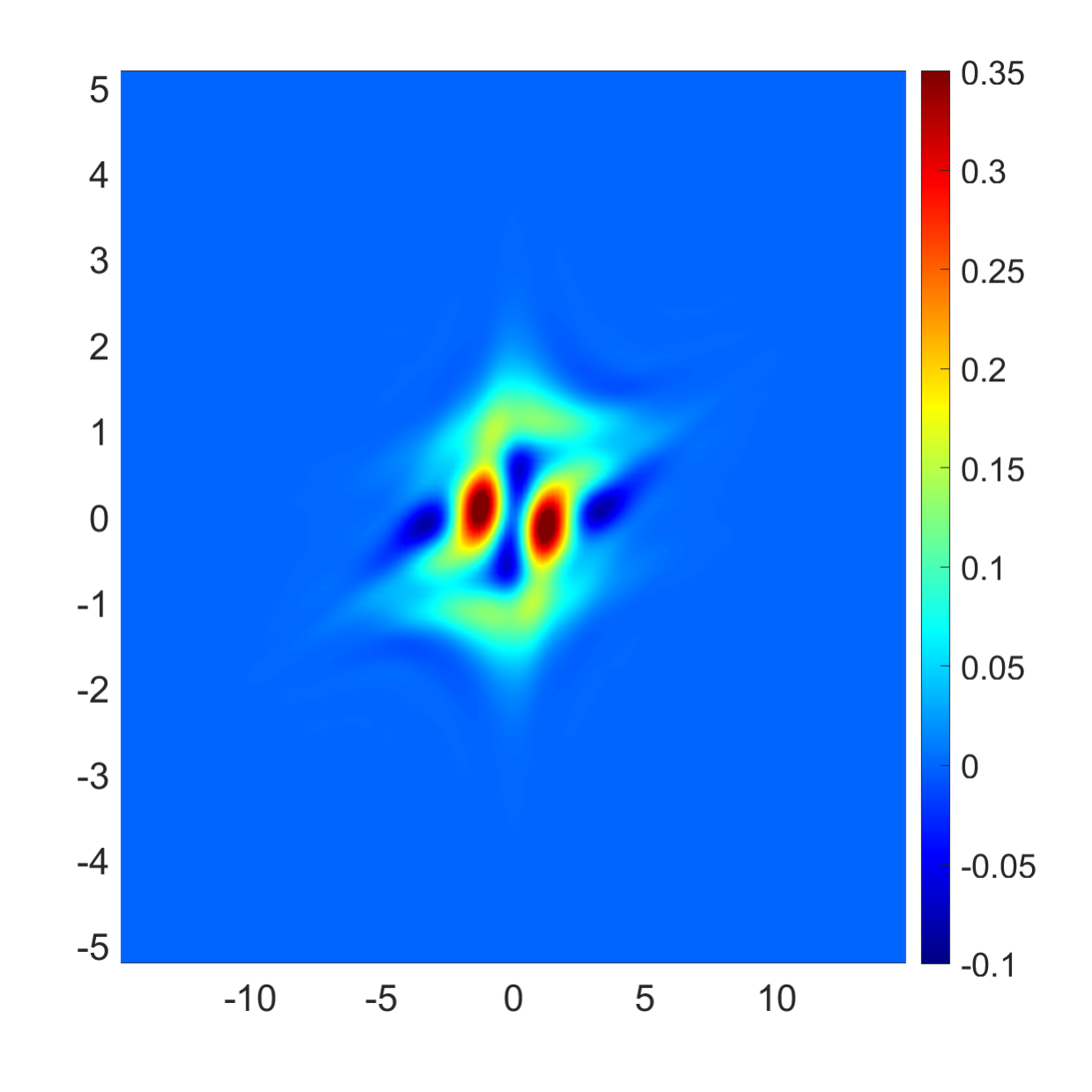}}
	\subfigure[t=6.]
	{\includegraphics[width=0.32\textwidth,height=0.25\textwidth]{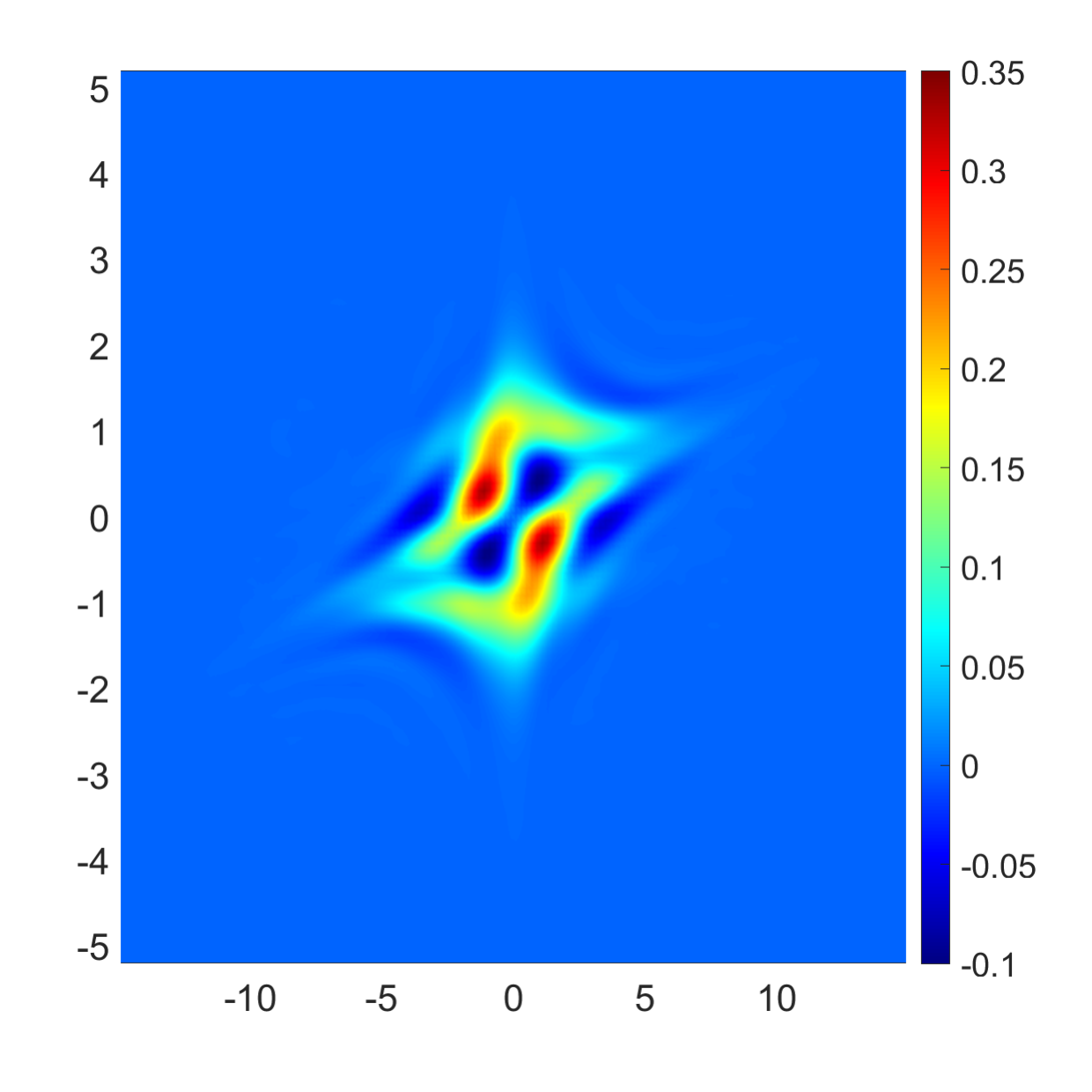}}\\
	\subfigure[t=7.]
	{\includegraphics[width=0.32\textwidth,height=0.25\textwidth]{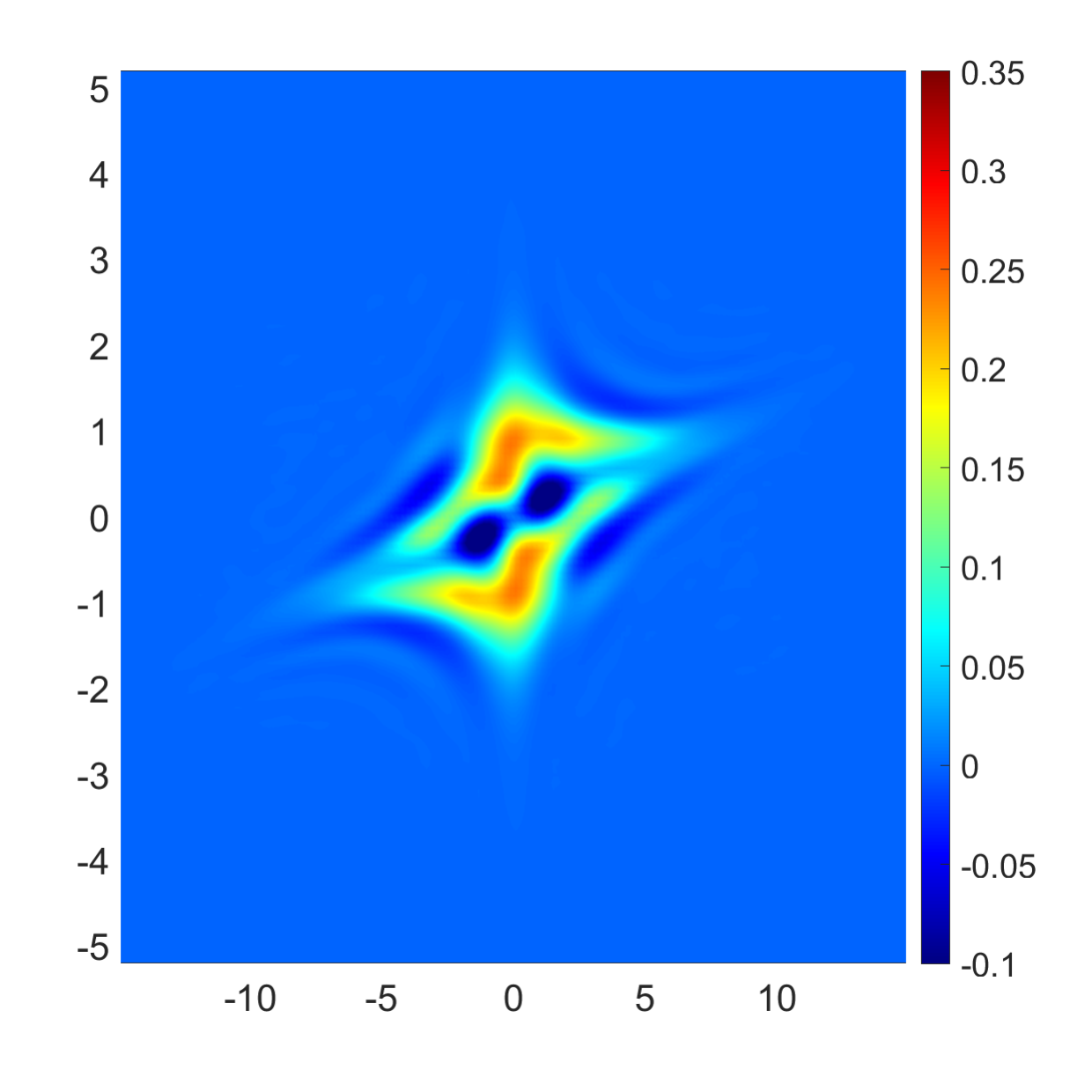}}
	\subfigure[t=8.]
	{\includegraphics[width=0.32\textwidth,height=0.25\textwidth]{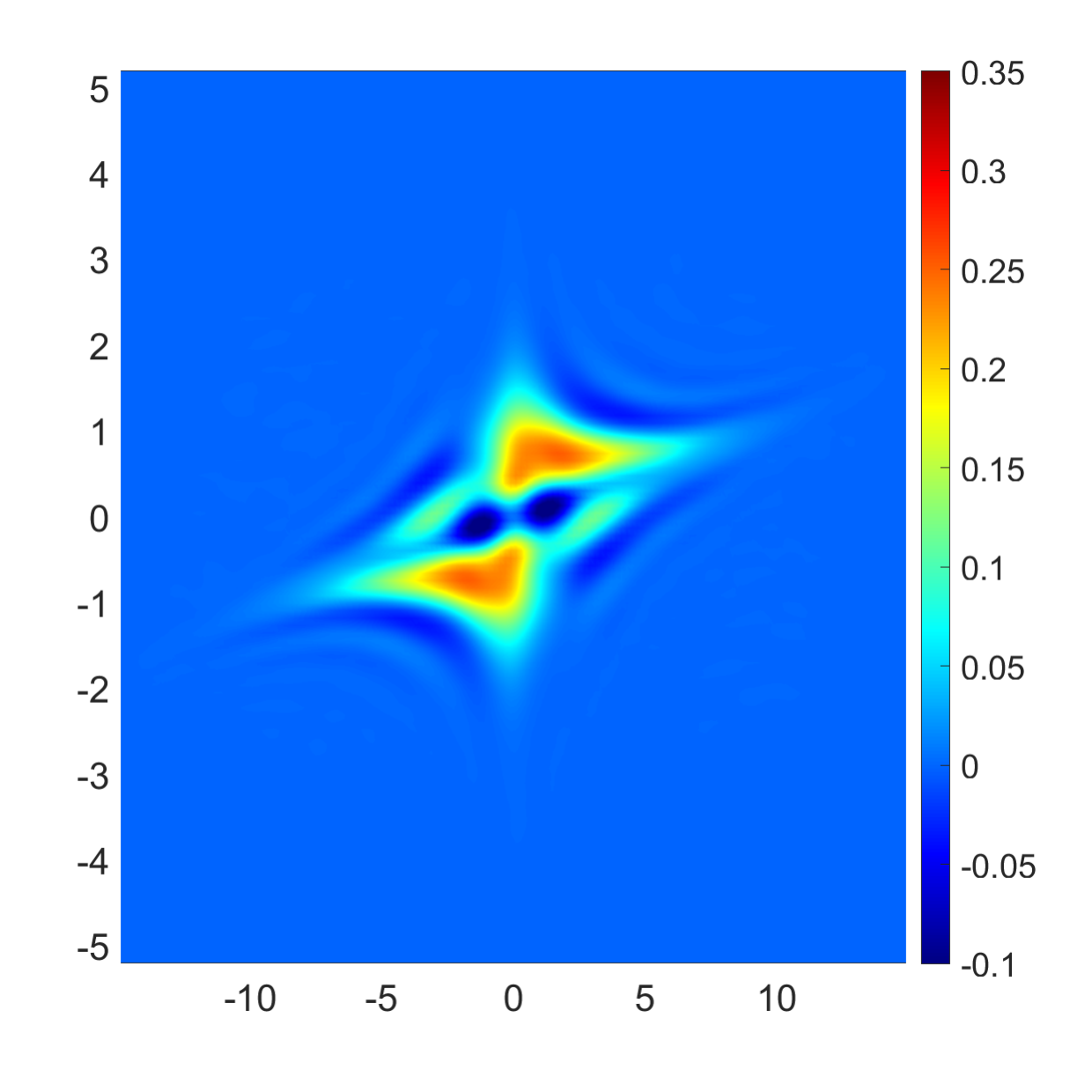}}
	\subfigure[t=9.]
	{\includegraphics[width=0.32\textwidth,height=0.25\textwidth]{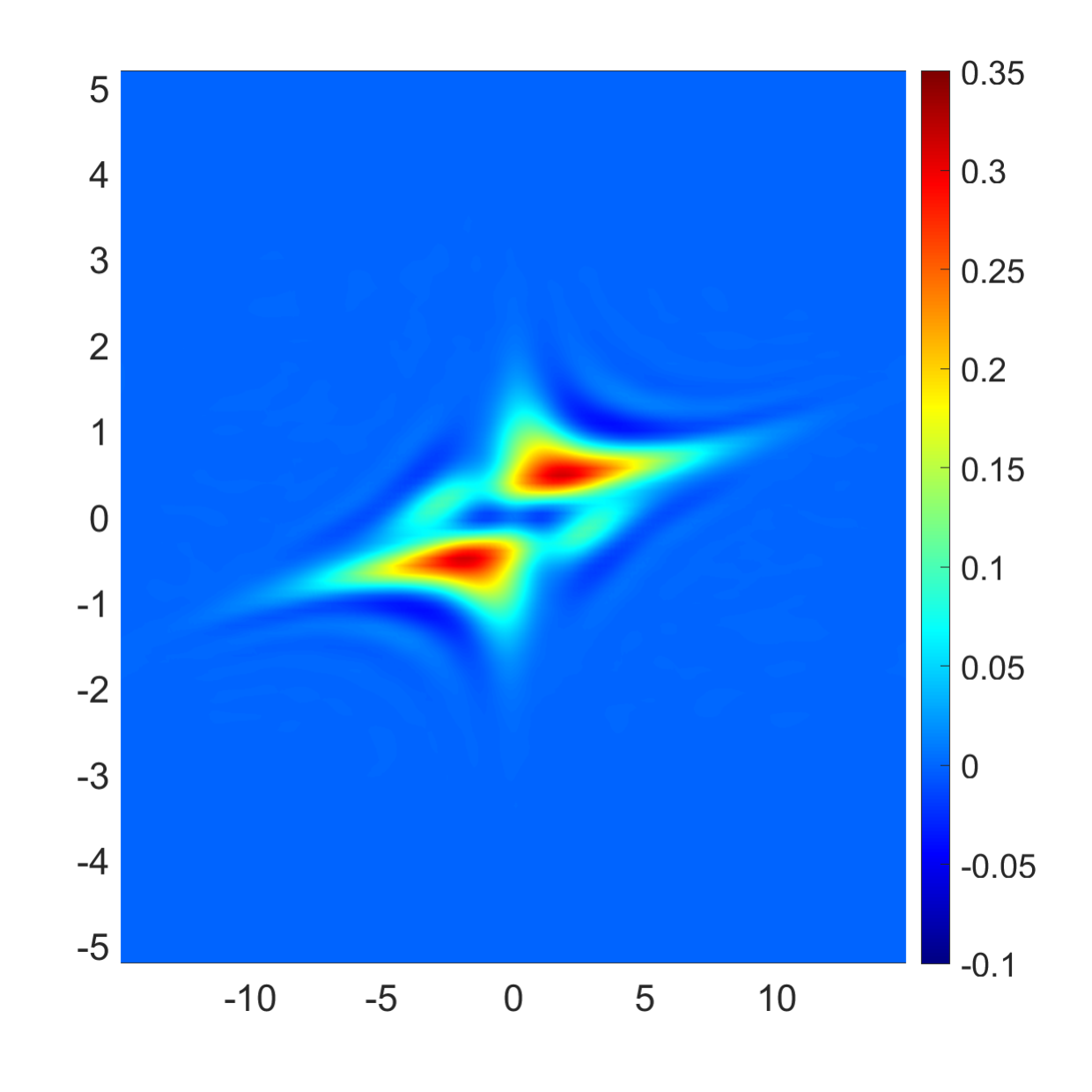}}\\
	\subfigure[t=10.]
	{\includegraphics[width=0.32\textwidth,height=0.25\textwidth]{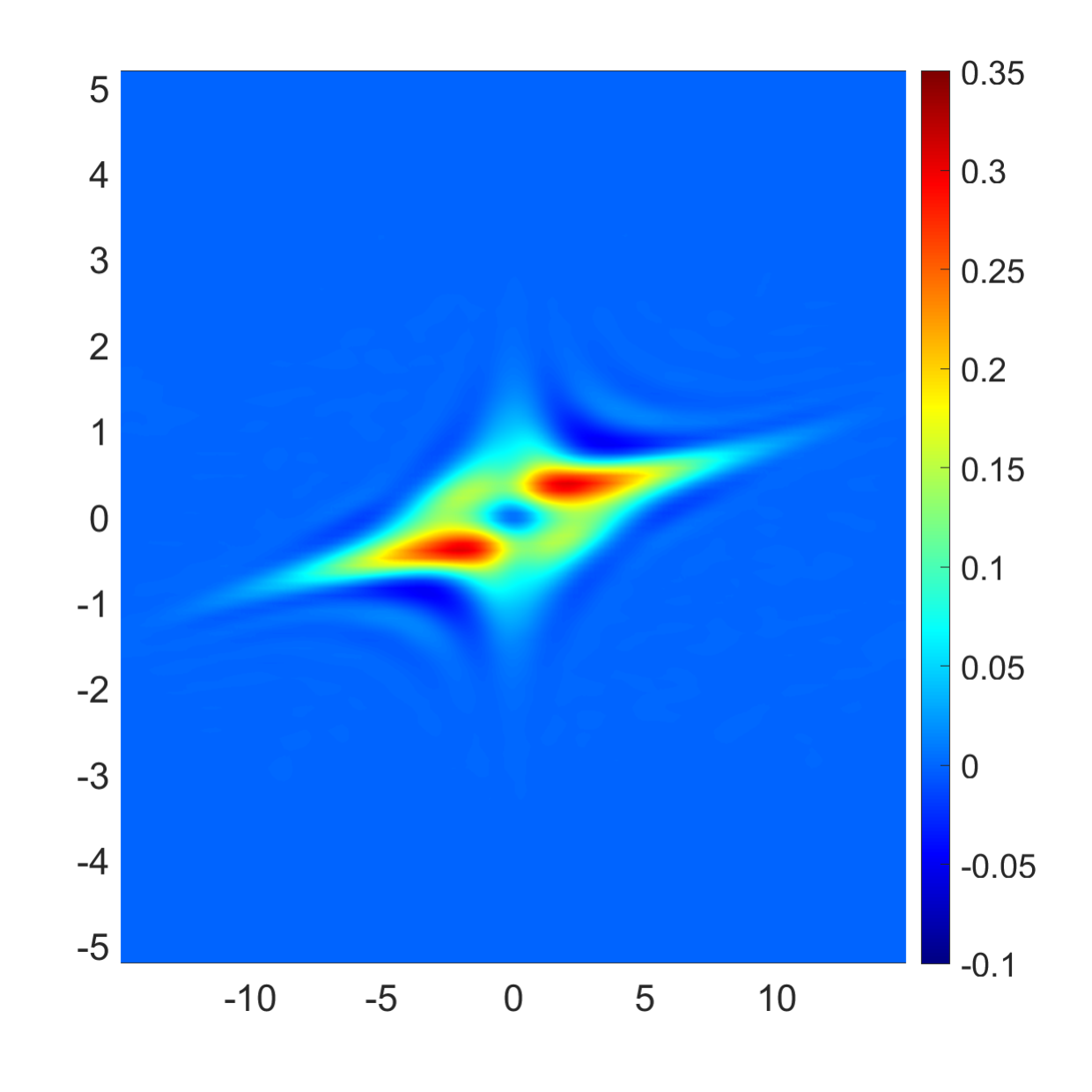}}
	\subfigure[t=11.]
	{\includegraphics[width=0.32\textwidth,height=0.25\textwidth]{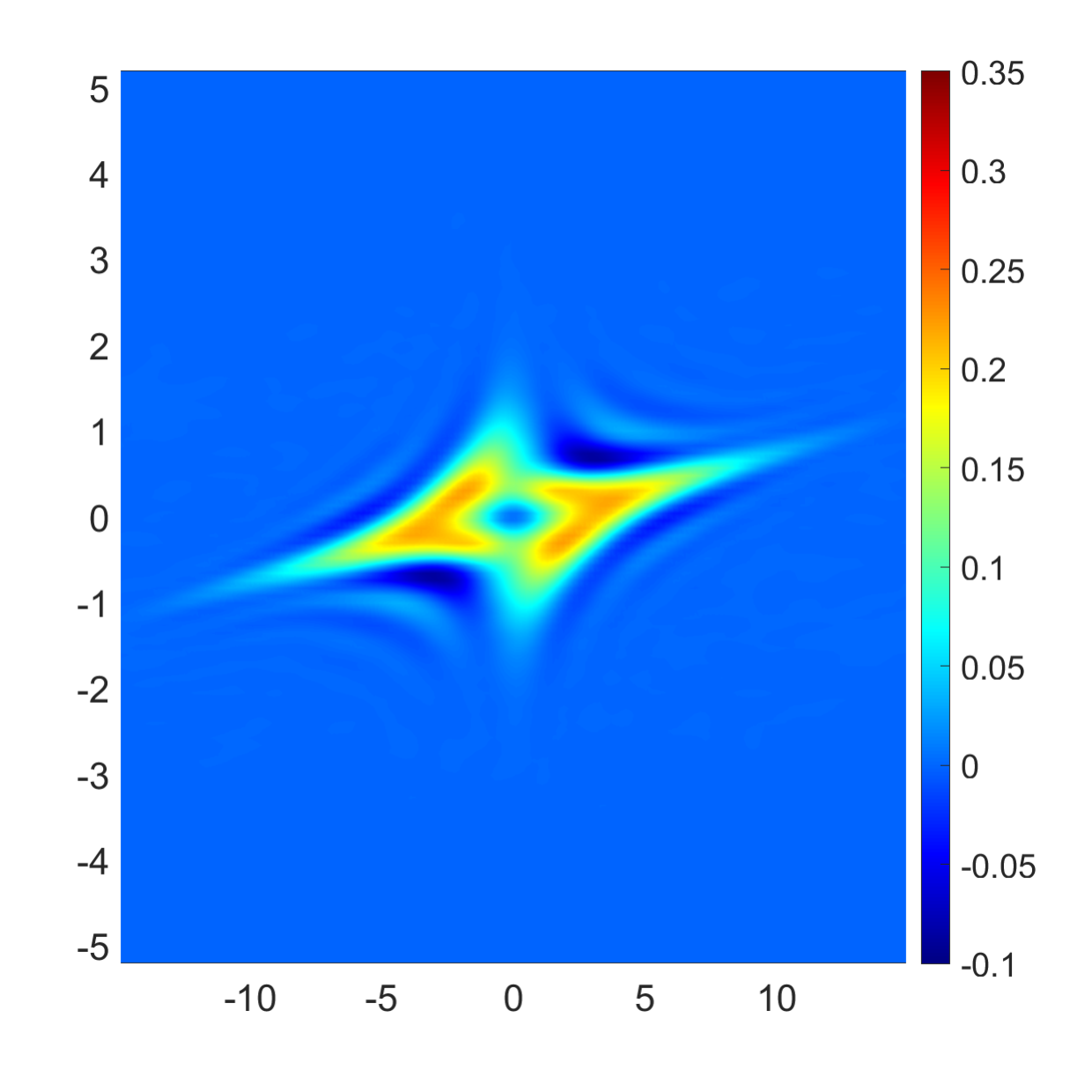}}
	\subfigure[t=12.]
	{\includegraphics[width=0.32\textwidth,height=0.25\textwidth]{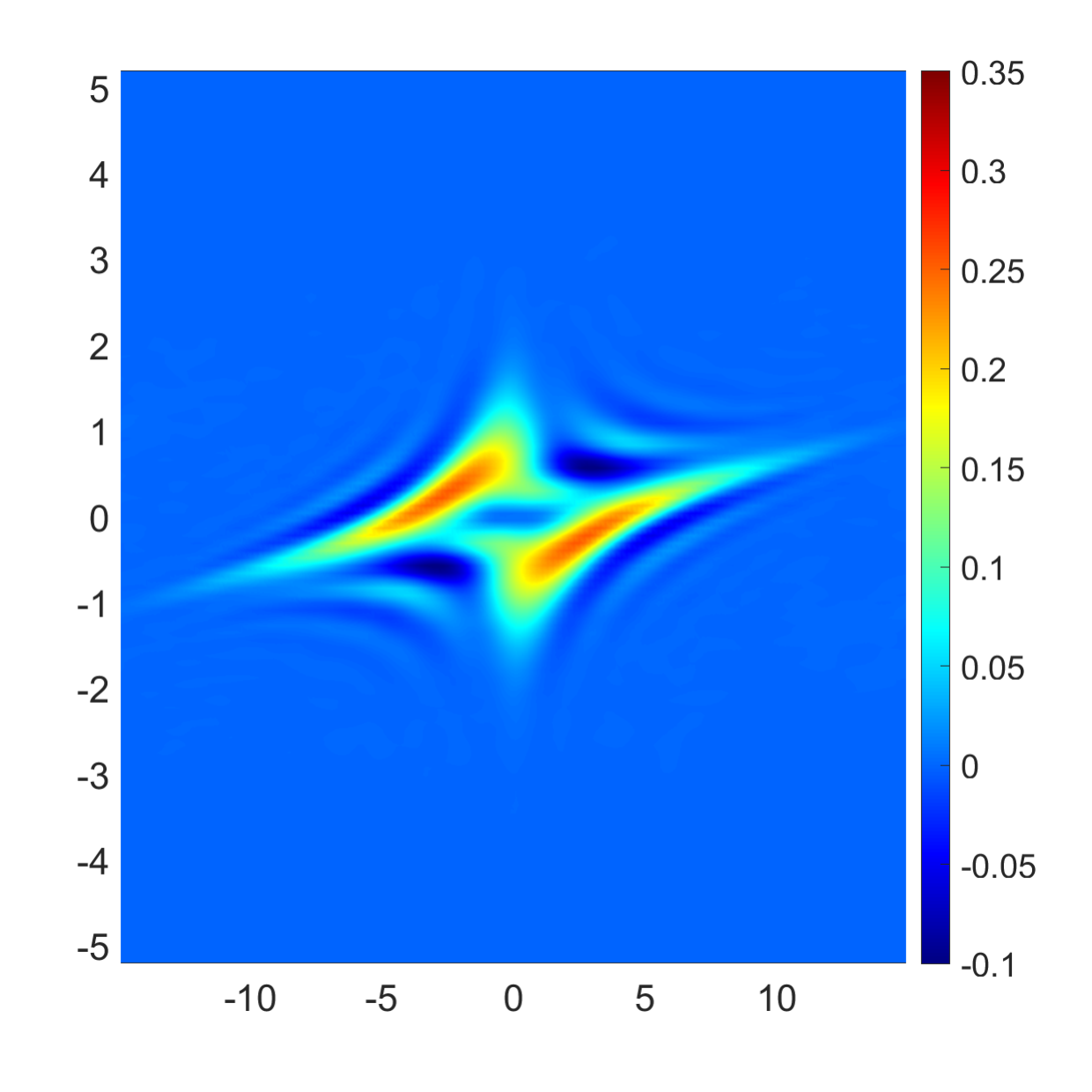}}\\
	\subfigure[t=13.]
	{\includegraphics[width=0.32\textwidth,height=0.25\textwidth]{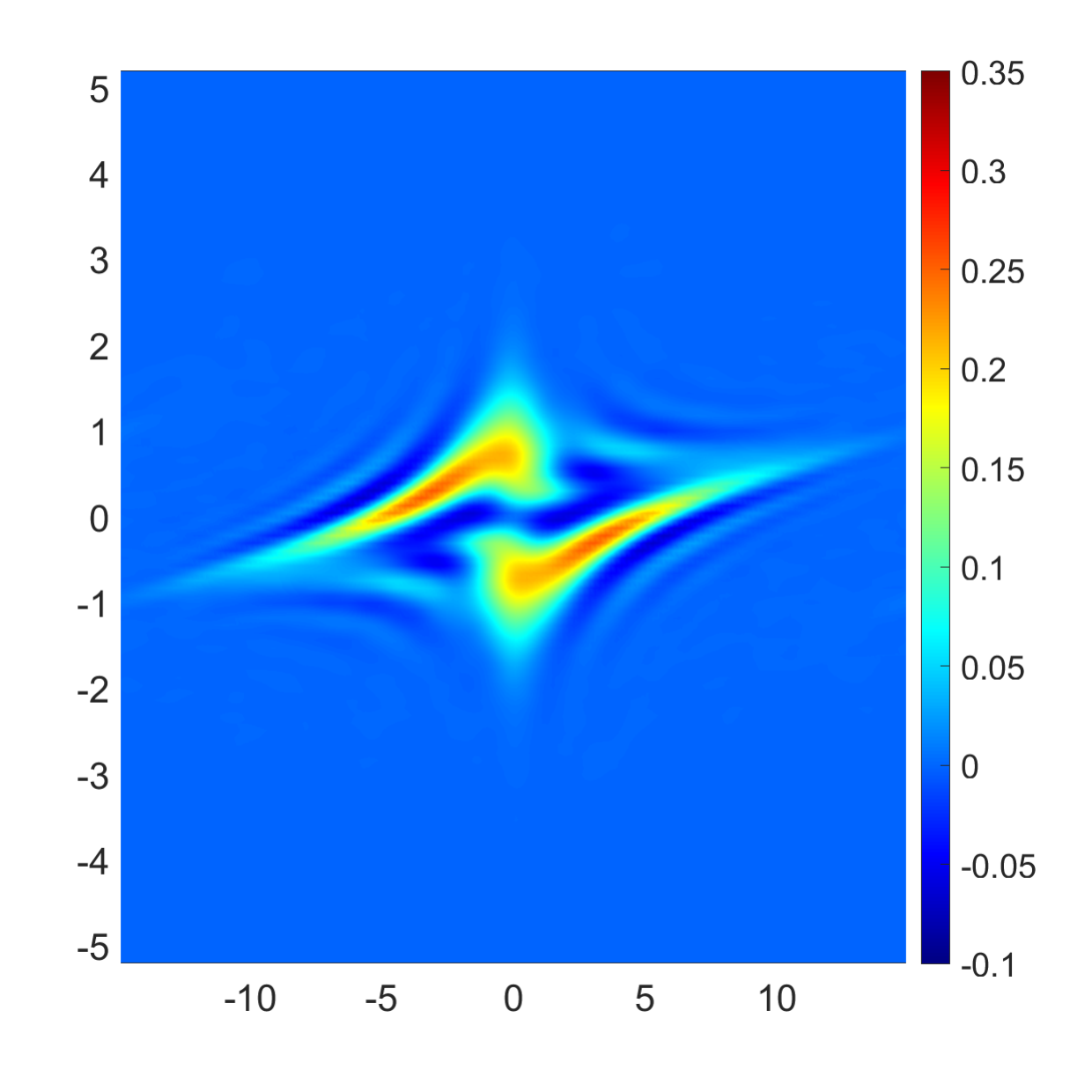}}
	\subfigure[t=14.]
	{\includegraphics[width=0.32\textwidth,height=0.25\textwidth]{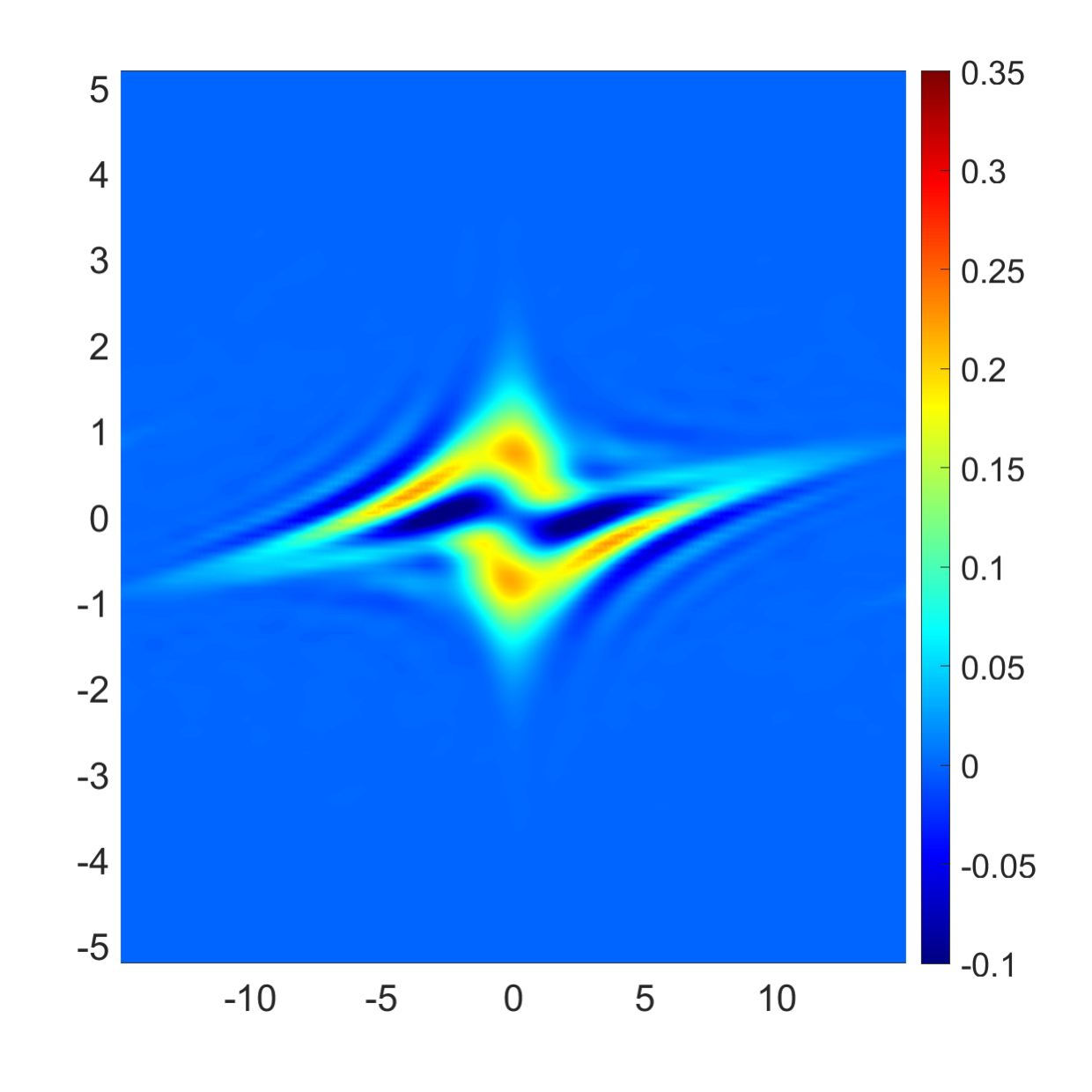}}
	\subfigure[t=15.]
	{\includegraphics[width=0.32\textwidth,height=0.25\textwidth]{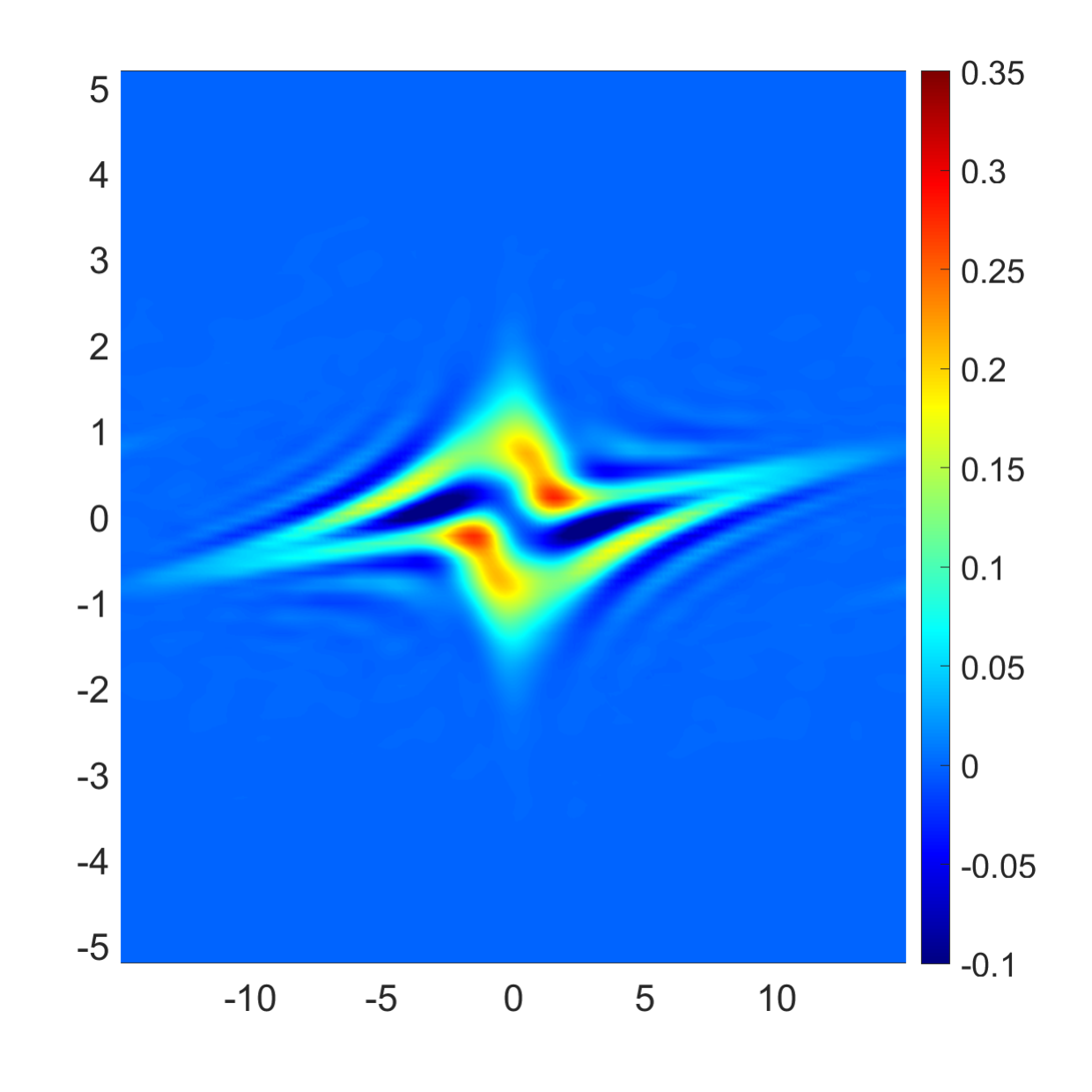}}
	\caption{The Helium-like system: $W_1(x,k,t)+W_2(x,k,t)$ at different times.}
	\label{fig::he}
\end{figure}

To demonstrate the dynamics of electrons more clearly, we take snapshots of the reduced Wigner function $W_1(x,k,t)+W_2(x,k,t)$ from $t=1$ to $t=15$, as shown \cref{fig::he}. 
At early times ($t<5$), the reduced Wigner function remains localized in the central region because of the nucleon–electron interaction. 
Afterwards, both dispersion and correlation become clearly visible, producing a highly oscillatory structure in phase space.
The region over which the reduced Wigner function evolves expands progressively with time, reaching the preset computational domain boundary by $t=15.$
It is worth noting that a region capable of taking negative values, referred to as the Fermi hole, persists around the centers of both wave packets. This is a direct consequence of the electron–electron interaction, since negative distributions are associated with areas inaccessible to experimental observation due to the constraints imposed by the uncertainty principle. 
The branching behavior of the Fermi hole varies at different times: for example, it exhibits two branches at $t=2,$ while only one branch is visible at $t=10.$ 
Nevertheless, it persists throughout, which agrees with the results from the ASM method \cite{xiong2016advectivespectralmixed} but differs from those obtained by the signed Monte Carlo method \cite{sellier2015simulation}.
This discrepancy may stem from the limited accuracy of the signed Monte Carlo method, which is insufficient to precisely capture quantum properties such as quantum interference.
This aspect warrants further investigation in the future.

\section{Conclusion and discussion}
\label{Conclusion and discussion}
This work presents the first application of DLRA to the Wigner equation. 
Our key innovation lies in introducing separable assumption \cref{DV low-rank} and combine it with two established approaches for $\pdo$ to transform $\pdo$ into a separated form \cref{DV alpha to theta}.
Comprehensive numerical tests on challenging non-stationary systems demonstrate both the accuracy and efficiency of our approach.
Even in the absence of a predetermined low-rank structure for the solution, DLRA can still serve as a numerical scheme that balances efficiency and accuracy.
These findings open new possibilities for simulating high-dimensional quantum systems.

While the current framework demonstrates promising results, several advanced strategies could further enhance the algorithm's performance. 
First, the second-order robust BUG integrator \cite{ceruti2024robust} provides stable second‑order accuracy and enables feasible adaptive rank control.
Second, changing the way of low-rank such as dynamical tensor approximation \cite{koch2010dynamical} and hierarchical low-rank approximation \cite{allmann-rahn2022parallel} may help reduce the required rank. 
Third, emerging approximate interpolation techniques \cite{chaturantabut2010nonlinear,christlieb2025samplingbased,dektor2025interpolatory} could be adopted to accelerate the computation.
Finally, beyond the extension in computational efficiency, we will consider more complex potentials in subsequent work, including singular potentials, discontinuous potentials, and even pure Coulomb potentials. 
These directions will be explored in our future research. 
\appendix


\end{document}